\documentclass[draft]{amsart}

\usepackage{amsmath,amssymb}
\usepackage{amsthm,amsrefs}
\usepackage{latexsym}
\usepackage{indentfirst}
\usepackage{mathrsfs}
\usepackage{graphicx}

\numberwithin{equation}{section}

\theoremstyle{plain}
\newtheorem{theorem}{Theorem}
\newtheorem*{theorem*}{Theorem}
\newtheorem{lemma}{Lemma}[section]

\newtheorem{prop}[lemma]{Proposition}
\newtheorem{coro}[lemma]{Corollary}

\theoremstyle{definition}
\newtheorem{defn}[lemma]{Definition}

\newtheorem{assump}{Assumption}

\theoremstyle{remark}
\newtheorem*{remark}{Remark}

\newcommand{\ud}{\,\mathrm{d}}
\newcommand{\RR}{\mathbb{R}}
\newcommand{\NN}{\mathbb{N}}
\newcommand{\ZZ}{\mathbb{Z}}

\newcommand{\LL}{\mathbb{L}}

\newcommand{\Or}{\mathcal{O}}

\newcommand{\wt}[1]{\widetilde{#1}}
\newcommand{\wh}[1]{\widehat{#1}}
\newcommand{\wb}[1]{\overline{#1}}

\newcommand{\mc}[1]{\mathcal{#1}}

\newcommand{\veps}{\varepsilon}

\newcommand{\x}{\times}

\newcommand{\abs}[1]{\lvert#1\rvert}
\newcommand{\norm}[1]{\lVert#1\rVert}

\newcommand{\set}[2]{\left\{\,#1\,\mid\,#2\,\right\}}

\newcommand{\inner}[2]{\left\langle#1,#2\right\rangle}
\newcommand{\Lr}[1]{\left(\,#1\,\right)}
\newcommand{\ffd}[1]{D_{#1}^+}
\newcommand{\bfd}[1]{D_{#1}^-}

\newcommand{\al}{\alpha}
\newcommand{\ga}{\gamma}
\newcommand{\pa}{\partial}
\newcommand{\na}{\nabla}
\newcommand{\dx}{\ud x}
\newcommand{\dt}{\ud t}

\newcommand{\wcb}{W_{\mathrm{CB}}}

\newcommand{\qc}{\mathrm{hy}}
\newcommand{\at}{\mathrm{at}}
\newcommand{\CB}{\mathrm{CB}}

\newcommand{\FE}{\mathrm{FE}}

\newcommand{\I}{\imath}
\newcommand{\sgn}{\operatorname{sgn}}
\renewcommand{\div}{\operatorname{\nabla\cdot}}

\newcommand{\barint}{\kern4pt \raise3.4pt\hbox{\vrule height.6pt
    width7pt} \kern-11pt \int}

\begin{document}

\title[Convergence of a force-based hybrid method]{Convergence of a
  force-based hybrid method for atomistic and continuum models in
  three dimension}

\author{Jianfeng Lu}
\address{Department of Mathematics \\
  Courant Institute of Mathematical Sciences \\
  New York University \\
  New York, NY 10012 \\
  email: jianfeng@cims.nyu.edu }

\author{Pingbing Ming}
\address{LSEC, Institute of Computational
  Mathematics and Scientific/Engineering Computing \\
  AMSS, Chinese Academy of Sciences \\
  No. 55, Zhong-Guan-Cun East Road \\
  Beijing 100190, China \\
  email: mpb@lsec.cc.ac.cn}

\date{February 8, 2011} \thanks{Part of the work was done during
  J.L.'s visit to State Key Laboratory of Scientific and Engineering
  Computing, Chinese Academy of Sciences. J.L. appreciates its
  hospitality. The work of P.B.M. was supported by National Natural
  Science Foundation of China under grants 10871197, 10932011, and by
  the funds from Creative Research Groups of China through grant
  11021101, and by the support of CAS National Center for Mathematics
  and Interdisciplinary Sciences. We thank Weinan E and Robert V. Kohn
  for helpful discussions.}

\begin{abstract}
  We study a force-based hybrid method that couples atomistic models
  with nonlinear Cauchy-Born elasticity models. We show that the
  proposed scheme converges quadratically to the solution of the
  atomistic model, as the ratio between lattice parameter and the
  characteristic length scale of the deformation tends to
  zero. Convergence is established for general short-ranged atomistic
  potential and for simple lattices in three dimension. The
  convergence is based on consistency and stability analysis. General
  tools are developed in the framework of pseudo-difference operators
  for stability analysis in arbitrary dimension of the multiscale
  atomistic and continuum coupling methods.
\end{abstract}

\maketitle

\section{Introduction}

Multiscale methods for mechanical deformation of materials have been
investigated intensely in recent years. The main spirit of these
methods is to use atomistic models for regions containing defects,
and continuum models in regions where the material is smoothly
deformed. We refer to the recent review~\cite{MillerTadmor:2009} for
various methods and the book \cite{E:book} for general discussion of
multiscale modeling.

There are two different ways of coupling atomistic and continuum
models. One is based on energy, and the other is based on force.
The energy-based method defines an energy which is a mixture of
atomistic energy and continuum elasticity energy. The energy
functional is then minimized to obtain the solution. The force-based
method works instead at the level of force balance equations. The
forces derived from atomistic and continuum models are coupled
together. The force balance equations are solved to obtain the
deformed state of the system.

From a numerical analysis point of view, one of the key issues for
these multiscale methods is the consistency and stability of the
coupled schemes. Taking one of the most successful
multiscale methods, the quasicontinuum
method~\cites{TadmorOrtizPhillips:1996, KnapOrtiz:2001} for example,
one of the main issues is the so called ghost force
problem~\cite{ShenoyMillerTadmorPhillipsOrtiz:1999}, which are the
artificial non-zero forces that the atoms experience at their
equilibrium state. In the language of numerical analysis, it means
that the scheme lacks consistency at the interface between atomistic
and continuum regions \cite{ELuYang:2006}. In~\cite{MingYang:2009},
it was shown that the ghost forces may lead to a finite size error
of the gradient of the solution.

The stability analysis for the coupling schemes is so far limited to
one dimensional systems, in which case a direct calculation is
possible thanks to the easy one dimensional lattice structure and
pairwise interaction potential. This is no longer the case in two and
three dimensions, and the extension is by on means easy. More general
tools for stability analysis are needed, to address in general the
multiscale hybrid methods.

In this work, based on existing ideas in the literature, we formulate
a force-based hybrid scheme for general short-ranged potentials (with
some natural assumptions) in three dimension. We focus on the
numerical analysis of the hybrid method, which is a representative of
a general class of multiscale methods. The solution of the proposed
method converges quadratically to the solution of the atomistic model
as the ratio between lattice parameter and the characteristic length
scale of the mechanical deformation goes to zero. To the best of our
knowledge, this is the first convergence result for multiscale methods
coupling atomistic and continuum models in three dimension.

The convergence result is based on the analysis of consistency and
linear stability. To achieve this, we study the linearized operator in
the framework of pseudo-difference operators. We obtained the
stability estimate combining regularity estimate of pseudo-difference
operators, consistency of the linearized operator, and stability of
the continuous problem. These tools developed will help understanding
multiscale methods in general.

Before we present the formulation of the method and the main theorem
in Section~\ref{sec:formulation}, we need some preliminaries and
notations.
\subsection{Lattice function and norms}

We will consider only Bravais lattices (see for
example~\cite{AshcroftMermin:1976}) in this work, denoted as $\LL$.
Let $d$ be the dimension. Let $\{a_j\} \subset \RR^d$, $j=1, \cdots,
d$ be basis vectors of the lattice $\LL$, hence
\[
  \LL = \{ x \in \RR^d \mid x = \sum_j n_j a_j,\, n \in \mathbb{Z}^d \}.
\]
Let $\{b_j\} \subset \RR^d$, $j = 1, \cdots, d$ be the reciprocal
basis vectors, given by
\[
a_j \cdot b_k = 2\pi \delta_{jk}.
\]
The reciprocal lattice $\LL^{\ast}$ is then
\[
\LL^{\ast} = \{ x \in \RR^d \mid x = \sum_j n_j b_j, \, n \in
\mathbb{Z}^d \}.
\]
Denote the unit cells of $\LL$ and $\LL^{\ast}$ as
$\Gamma$ and $\Gamma^{\ast}$ respectively.
\begin{align*}
  & \Gamma = \{ x \in \RR^d \mid x = \sum_j c_j a_j,\; 0 \leq c_j < 1,
  \, j = 1, \cdots, d\}; \\
  & \Gamma^{\ast} = \{ x \in \RR^d \mid x = \sum_j c_j b_j,\; -1/2
  \leq c_j < 1/2,  \, j = 1, \cdots, d\}.
\end{align*}

For $\veps = 1/n, \, n \in \ZZ_+$, we will consider lattice system
$\veps \LL$ inside domain $\Omega = \Gamma \subset \RR^d$, denoted as
$ \Omega_{\veps} = \Omega \cap \veps \LL$.  Note that the lattice
constant is $\veps$, so that the number of points in $\Omega_{\veps}$
is $1/ \veps^d$. We will restrict to periodic boundary conditions in
this work, general boundary conditions will be leaved for future
publications. For a lattice function $u$ defined on $\veps \LL$, we
say it is \emph{$\Omega_{\veps}$-periodic} if
\[
  u(x) = u(x'), \qquad \forall \, x, x'
  \in \veps \LL,\, x - x' = a_j \text{ for some } j \in \{1, \cdots, d\}.
\]
In particular, an $\Omega_{\veps}$-periodic function is determined
by its restriction on $\Omega_{\veps}$. Functions defined on
$\Omega_{\veps}$ can be easily extended to $\Omega_{\veps}$-periodic
functions defined on $\veps\LL$.

We also define the reciprocal lattice associated with
$\Omega_{\veps}$. Let $\LL^{\ast}_{\veps} = \LL^{\ast} \cap (
\Gamma^{\ast} / \veps ) $. Define $K_{\veps}$ a subset of $\ZZ^d$
given by
\[
  K_{\veps} = \{ \mu \in \ZZ^d \mid \sum_j \veps \mu_j b_j \in \Gamma^{\ast} \},
\]
hence $\LL^{\ast}_{\veps}$ is given by
\[
  \LL^{\ast}_{\veps} = \{ x \in \RR^d \mid x = \sum_j \mu_j b_j,\,
  \mu \in K_{\veps}\}.
\]

For $\mu \in \ZZ^d$, the translation operator $T^{\mu}_{\veps}$ is
defined as
\[
  (T^{\mu}_{\veps} u)(x) = u(x + \veps \mu_j a_j), \quad \text{for } x \in \RR^d.
\]
We define the forward and backward discrete gradient operators as
\[
\ffd{\veps, s} =\veps^{-1}(T^{\mu}_{\veps}-I)\qquad\text{and}\qquad
\bfd{\veps, s}=\veps^{-1}(I-T^{\mu}_{\veps}),
\]
where $s = \sum_i \mu_i a_i$ and $I$ denotes the identity operator. It
is easy to see $\ffd{\veps, -s}=-\bfd{\veps, s}$.

We say $\alpha$ is a multi-index, if $\alpha \in \ZZ^d$ and $\alpha
\geq 0$. We will use the notation
\[
  \abs{\alpha} = \sum_{j=1}^d \alpha_j.
\]
For a multi-index $\alpha$, the difference operator
$D_{\veps}^{\alpha}$ is given by
\[
D^{\alpha}_{\veps} = \prod_{j=1}^d (\ffd{\veps, a_j})^{\alpha_j}.
\]
When no confusion will occur, we will omit the subscript $\veps$ in
the notations $T^{\mu}_{\veps}$, $\ffd{\veps, s}$, $\bfd{\veps, s}$
and $D^{\alpha}_{\veps}$ for simplicity.

We will use various norms for functions defined on the lattice
$\Omega_{\veps}$. For integer $k\geq 0$, define the difference norm
\[
  \norm{u}_{\veps, k}^2 = \sum_{0 \leq \abs{\alpha} \leq k} \veps^d
  \sum_{x \in \Omega_{\veps}} \abs{(D_{\veps}^{\alpha} u)(x)}^2.
\]
It is clear that $\norm{\cdot}_{\veps, k}$ is a discrete analog of
Sobolev norm associated with $H^k(\Omega)$. Hence, we denote the
corresponding spaces of lattice functions as $H_{\veps}^k(\Omega)$ and
$L_{\veps}^2(\Omega)$ when $k = 0$.  We also need the uniform norms on
the lattice $\Omega_{\veps}$, given by
\begin{align*}
  & \norm{u}_{L^{\infty}_{\veps}} = \max_{x \in \Omega_{\veps}}
  \abs{u(x)}, \\
  & \norm{u}_{W^{k,\infty}_{\veps}} = \sum_{0 \leq \abs{\alpha} \leq
    k} \max_{x \in \Omega_{\veps}} \abs{(D^{\alpha}_{\veps} u)(x)}.
\end{align*}
In the above definitions, we have identified lattice function $u$ with
its $\Omega_{\veps}$-periodic extension to function defined on $\veps
\LL$, and hence the differences are well-defined.  These norms extend
to vector-valued functions as usual.

Define the discrete Fourier transform for lattice functions $f$ as
\begin{equation}\label{eq:discF}
  \wh{f}(\xi) = \veps^d (2\pi)^{-d/2}
  \sum_{x \in \Omega_{\veps}} e^{-\I \xi \cdot x} f(x), \quad \xi \in
  \mathbb{L}^{\ast}_{\veps},
\end{equation}
and its inverse as
\begin{equation}\label{eq:invdiscF}
  f(x) = (2\pi)^{d/2} \sum_{\xi \in \mathbb{L}^{\ast}_{\veps}} e^{\I x \cdot \xi}
  \wh{f}(\xi), \quad x \in \Omega_{\veps}.
\end{equation}

We need a symbol which plays the same role for difference operators
that $\Lambda^2(\xi) = 1 + \Lambda_0^2(\xi) = 1 + \abs{\xi}^2$ plays
for differential operators. For $\veps > 0, \, \xi \in
\LL_{\veps}^{\ast}$, let
\[
  \Lambda_{j, \veps}(\xi) = \frac{1}{\veps} \abs{ e^{\I\veps
      \xi_j} - 1}, \qquad j = 1, \cdots, d,
\]
and
\[
\Lambda_{\veps}^2(\xi) = 1 + \Lambda_{0, \veps}^2(\xi) = 1 +
\sum_{j=1}^d \Lambda_{j, \veps}^2(\xi) = 1 + \sum_{j=1}^d
\frac{4}{\veps^2} \sin^2\Bigl( \frac{\veps \xi_j}{2} \Bigr).
\]
It is not hard to check for any $\xi \in \LL_{\veps}^{\ast}$, it holds
\begin{equation}\label{eq:symbolcompare}
  c\Lambda^2(\xi)
  \leq \Lambda_{\veps}^2(\xi) \leq \Lambda^2(\xi).
\end{equation}
where the positive constant $c$ depends on $\{b_j\}$.

The $L^2_{\veps}$ norm of lattice function can be rewritten as
\begin{equation}
  \norm{f}_{\veps, 0}^2 = (2\pi)^d
  \sum_{\xi \in \LL_{\veps}^{\ast}} \abs{\wh{f}(\xi)}^2.
\end{equation}
Indeed, using Poisson summation formula,
\begin{equation*}
  \begin{aligned}
    \sum_{\xi \in \LL_{\veps}^{\ast}} \abs{\wh{f}(\xi)}^2 & =
    \sum_{\xi\in\LL_{\veps}^{\ast}} \veps^{2d} (2\pi)^{-d} \sum_{x \in
      \Omega_{\veps}} e^{\I \xi\cdot x} f^{\ast}(x) \sum_{x' \in
      \Omega_{\veps}} e^{-\I \xi \cdot x'}f(x') \\
    & = \sum_{x, x' \in \Omega_{\veps}} \veps^{2d} (2\pi)^{-d}
    f^{\ast}(x) f(x') \sum_{\xi \in \LL_{\veps}^{\ast}} e^{\I \xi \cdot(x - x')} \\
    & = \sum_{x \in \Omega_{\veps}} (2\pi)^{-d} \veps^d \abs{f(x)}^2 =
    (2\pi)^{-d}\norm{f}_{\veps, 0}^2.
  \end{aligned}
\end{equation*}
Moreover, notice that for $\xi \in \LL_{\veps}^{\ast}$, we have
\[
  \wh{D_{\veps, a_j}^+ f}(\xi) = \frac{1}{\veps} ( e^{\I \veps \xi \cdot a_j} - 1)
  \wh{f}(\xi).
\]
Therefore, discrete Sobolev norms have equivalent representations
using discrete Fourier transform:
\[
  c \norm{f}_{\veps, k}^2 \leq \sum_{\xi \in \LL_{\veps}^{\ast}} \Lambda_{\veps}^k(\xi)
  \abs{\wh{f}(\xi)}^2 \leq  C \norm{f}_{\veps, k}^2,
\]
with positive constant $c$ depending on $k$ and $\{ a_j \}$.

For $k>d/2$, we have the following discrete Sobolev imbedding
inequality~\cite{Frank:1971}*{Proposition 6}:
\[
\norm{f}_{L^{\infty}_{\veps}}\le C\norm{f}_{\veps,k},
\]
where $C$ depends on $k$ and $\Omega$.
\subsection{Atomistic model and Cauchy-Born rule}

In this work, we will restrict our attention to classical empirical
potentials.  For atoms located at $\{ y_1, \cdots, y_N\}$, the
interaction potential energy between the atoms is given by
\[
  V(y_1,\cdots,y_N),
\]
where $V$ often takes the form:
\[
V(y_1,\cdots,y_N)=\sum_{i,j}V_2(y_i/\veps, y_j/\veps)
+\sum_{i,j,k}V_3(y_i/\veps,y_j/\veps,y_k/\veps)+\cdots.
\]
Here we have omitted interactions of more than three atoms.

Different potentials are chosen for different materials. In this
paper, we will work with general atomistic models, and we will make
the following assumptions on the potential functions $V$ as
in~\cite{EMing:2007}:
\begin {enumerate}

\item $V$ is translation invariant.

\item $V$ is invariant with respect to rigid body motion.

\item $V$ is smooth in a neighborhood of the equilibrium state.

\item $V$ has finite range and consequently we will consider only
interactions that involve a finite number of atoms.
\end {enumerate}
The first two assumptions are general~\cite{BornHuang:1954}, while
the latter two are specific technical assumptions.

In fact, for simplicity of notation and clarity of presentation, our
presentation will be limited to potentials that contain only two-body
and three-body potentials. Actually, we will sometimes only make
explicit the three-body terms in the expressions for the potential and
omit the two-body terms. It is straightforward to extend our results
to potentials with interactions of more atoms that satisfy the above
conditions, following the discussion on the three-body terms.
By~\cite{Keat:1965}, the potential function $V$ is a function of atom
distances and angles by invariance with respect to rigid body
motion. Therefore, we may write
\begin{align*}
  V_2(y_i,y_j)&=V_2\Lr{\abs{y_i-y_j}^2},\\
  V_3(y_i,y_j,y_k)& =V_3\Lr{\abs{y_i-y_j}^2,
    \abs{y_i-y_k}^2,\inner{y_i-y_j}{y_i-y_k}},
\end{align*}
where $\inner{\cdot}{\cdot}$ denotes the inner product over $\RR^d$.
We write the two-body and three-body potentials this way to make the
formula in our calculations easier to read.

We assume that the atoms are located at $\Omega_{\veps}$ in
equilibrium, with $x$ denoting the equilibrium position ($x \in
\Omega_{\veps}$). Positions of the atoms under deformation will be
viewed as a function defined over $\Omega_{\veps}$, denote as $y(x) =
x + u(x)$. Hence, $u: \Omega_{\veps} \to \RR^d$ is the displacement of
atoms. We extend $u$ as an $\Omega_{\veps}$-periodic function defined
on $\veps \LL$. Denote the space of atom positions $y$ as
\begin{equation*}
  X_{\veps} = \{ y: \veps \LL \to \RR^d
  \mid y = x + u,\, u\ \Omega_{\veps}\text{-periodic},\, 
  \sum_{x\in\Omega_{\veps} u(x) = 0} \}.
\end{equation*}
Hence, $y \in X_{\veps}$ satisfies
\begin{equation*}
  y(x) - y(x') = x - x', \qquad \forall \, x, x'
  \in \veps \LL,\, x - x' = a_j \text{ for some } j \in \{1, \cdots, d\}.
\end{equation*}

The atomistic problem is formulated as follows. For given $f:
\Omega_{\veps} \to \RR^d$, find $y \in X_{\veps}$ such that
\begin{equation}\label {atom:min}
  y = \arg\min_{z \in X_{\veps}} I_{\at}(z),
\end{equation}
where
\begin{equation*}
  I_{\at}(z) = \dfrac{1}{3!}  \veps^d \sum_{x\in\Omega_{\veps}}\sum_{(s_1,s_2)
    \in S} V_{(s_1, s_2)}[z] - \veps^d \sum_{x \in
    \Omega_{\veps}} f(x) z(x),
\end{equation*}
where
\begin{equation*}
  V_{(s_1, s_2)}[z] =
  V\Lr{\abs{\ffd{s_1}z(x)}^2,
    \abs{\ffd{s_2}z(x)}^2,\inner{\ffd{s_1}z(x)}{\ffd{s_2}z(x)}}.
\end{equation*}
Here $S$ is the set of all possible $(s_1, s_2)$ within the range of
the potential. By our assumptions, $S$ is a finite set.  Note that as
remarked above, we only make explicit the three body terms in the
potential. In $I_{\at}$, $\veps^d$ is a normalization factor, so that
$I_{\at}$ is actually the energy of the system per atom.

The Euler-Lagrange equations for the atomistic problem is then
\begin{equation}\label{atom:eq}
  \mc{F}_{\at}[y](x) = f(x),\qquad x\in\Omega_{\veps},
\end{equation}
where
\begin{align*}
    \mc{F}_{\at}[y](x) &= \sum_{(s_1,s_2)\in S}
    \Bigl( \bfd{s_1}\Lr{2\pa_1
      V_{(s_1, s_2)}[y](x)\ffd{s_1}y(x)+\pa_3V_{(s_1, s_2)}[y](x)\ffd{s_2}y(x)}\\
    &\phantom{\sum_{(s_1,s_2)}}\qquad+\bfd{s_2}\Lr{2\pa_2 V_{(s_1,
        s_2)}[y](x)\ffd{s_2}y(x)+\pa_3V_{(s_1,
        s_2)}[y](x)\ffd{s_1}y(x)}\Bigr),
\end{align*}
where for $i=1,2,3$, we denote
\[
\pa_i V_{(s_1, s_2)}[y](x)=\partial_i
V\Lr{\abs{\ffd{s_1}y(x)}^2,\abs{\ffd{s_2}y(x)}^2,
  \inner{\ffd{s_1}y(x)}{\ffd{s_2}y(x)}},
\]
the partial derivative with respect to the $i$-th argument of $V$.

To introduce the continuum Cauchy-Born (CB) elasticity problem
~\cites{BornHuang:1954, Ericksen:1984, Ericksen:2008}, we fix more
notations. For any positive integer $k$, we denote by
$W^{k,p}(\Omega;\RR^d)$ the Sobolev space of mappings
$y{:}\;\Omega\to\RR^d$ such that $\|y\|_{W^{k,p}}<\infty$. In
particular, $W_{\sharp}^{k,p}(\Omega;\RR^d)$ denotes the Sobolev
space of periodic functions whose distributional derivatives of
order less than $k$ are in the space $L^p(\Omega)$.  For any $p>d$
and $m\ge 0$, we define $X$ as
\[
X=\{y: \Omega \to \RR^d \mid y = x + v,\, v\in
W^{m+2,p}(\Omega;\RR^d) \cap\,W_{\sharp}^{1,p}(\Omega;\RR^d), \,
\int_{\Omega}v=0 \}.
\]

As in~\cite{EMing: 2007}, we have the Cauchy-Born elasticity
problem as: find $y\in X$ such that
\begin{equation}\label {cb:min}
  y = \arg \min_{z\in X}I(z),
\end{equation}
where the total energy functional $I$ is given by
\[
I(z)=\int_{\Omega} \Lr{\wcb(\na v(x))-f(x)z(x)} \dx,
\]
where $v(x) = z(x) - x$ and Cauchy-Born stored energy density $\wcb$
is given by
\[
\wcb(A)=\dfrac{1}{3!}\sum_{(s_1,s_2)\in S} W_{(s_1, s_2)}(A),
\]
where for $A \in \RR^{d\times d}$,
\[
W_{(s_1, s_2)}(A) = V\Lr{\abs{s_1+s_1
    A}^2,\abs{s_2+s_2A}^2,\inner{s_1+s_1A}{s_2+s_2A}}.
\]
The range $S$ is the same as that in the atomistic potential.  We have
used the deformed position $y$ instead of the more usual displacement
field $u$ as variable in \eqref{cb:min} in order to be parallel with
the atomistic problem.

The Euler-Lagrange equation for the Cauchy-Born elasticity
model is then
\begin {equation}\label {cb:eq}
  \mc{F}_{\CB}[y](x)=f(x),
\end {equation}
where
\begin{equation*}
  \mc{F}_{\CB}[y](x)=-\div\Lr{D_A\wcb(\na v(x))},
  \qquad v(x)=y(x)-x.
\end{equation*}
Here $D_A \wcb(A)$ denotes differentiation of $\wcb(A)$ with respect
to $A$.

Since we are primarily interested in the coupling between the
atomistic and continuum region, we will take the finite element
discretization $\mathcal{T}_{\veps}$ be a triangulation of
$\Omega_{\veps}$ with each atom site as an element vertex with element
size $\veps$. The triangulation is chosen so that it is translation
invariant. The approximation space $\wt{X}_{\veps}$ is defined as
\[
\wt{X}_{\veps}=\bigl\{ y \in W^{1,p}_{\sharp}(\Omega;\RR^d) \mid
y|_T\in P_1(T), \ \forall\, T\in\mc{T}_{\veps}\bigr\},
\]
where $P_1(T)$ is the space of linear functions on the element $T$.
\subsection{Force-based hybrid method}\label{sec:formulation}

We are ready to formulate the force-based hybrid method.

We take $\varrho: \Omega \to [0, 1]$ as a smooth standard cutoff
function. The atomistic region corresponds to the zero level set of
$\varrho$: $\Omega_{a} = \{x \mid \varrho(x) = 0\}$, and the
continuum region corresponds to the region that $\varrho$ equals to
$1$: $\Omega_{c} = \{ x \mid \varrho(x) = 1 \}$. The region in
between is a buffer between the atomistic and continuum regions.

The force-based hybrid method is given as: find $y(x)\in
X_{\veps}$ such that
\begin{equation}\label {sqc:eq}
  \mc{F}_{\qc}[y](x)\equiv (1 - \varrho(x)) \mc{F}_{\at}[y](x)
  + \varrho(x) \mc{F}_{\veps}[y](x) = f(x),\qquad
  x\in\Omega_{\veps},
\end{equation}
where $\mc{F}_{\veps}$ is the force from finite element
approximation of Cauchy-Born elasticity problem~\eqref{cb:min}. Due
to the choice of $\varrho$, in the atomistic region $\Omega_{a}$,
the force acting on the atom is just that of atomistic model, while
in the continuum region $\Omega_c$, the force is calculated from
finite element approximation of the Cauchy-Born elasticity.

The proposed scheme works in dimension $d \leq 3$ for general
short-range interaction potentials. The main result for this work is
the following quadratic convergence result for the force-based hybrid
method.
\begin{theorem}[Convergence]\label{thm:main}
  Under Assumptions~\ref{assump:stabatom} and \ref{assump:stabCB},
  there exist positive constants $\delta$ and $M$, so that for any $p
  > d$ and $ f \in W^{15, p}(\Omega) \cap W^{1, p}_{\sharp}(\Omega) $
  with $\norm{f}_{W^{15, p}} \leq \delta$, we have
  \begin{equation}\label{eq:main}
    \norm{y_{\qc} - y_{\at}}_{\veps,2} \leq M \veps^2.
  \end{equation}
\end{theorem}

\begin{remark}
  While we do not attempt in this work to optimize the regularity
  assumption on $f$, we note that it is easy to relax the assumption
  to $f \in W^{5, p}(\Omega)$ with $p>d$ following the remarks below
  in the proof.
\end{remark}

\begin{remark}
  The sharp stability conditions Assumptions~\ref{assump:stabatom} and
  \ref{assump:stabCB} will be given in
  Section~\ref{sec:regularity}. These assumptions are quite natural
  and physical. We refer to Section~\ref{sec:regularity} and also
  \cite{EMing:2007} for more discussions on the stability conditions
  and its link to physics literature.
\end{remark}

The proof of Theorem~\ref{thm:main}, which will be viewed as a
convergence result for (nonlinear) finite difference schemes,
follows the spirit of Strang's work \cite{Strang:1964}. In short,
consistency and linear stability implies convergence. The heart of
the matter lies in the analysis of consistency and stability, which
will be the focus of the proof.

The rest of the paper is organized as follows. In the next
subsection, we review some related works.
Section~\ref{sec:consistency} discusses the consistency of the
scheme. The linear stability is proved in
Section~\ref{sec:stability}.  The stability estimate is based on the
regularity estimate of finite difference schemes in
Section~\ref{sec:regularity}, which is established in the framework
of pseudo-difference operators \cites{LaxNirenberg:1966,
Thomee:1964,BubeStrikwerda:1983}. With the preparation of
consistency and linear stability analysis, the proof is concluded in
Section~\ref{sec:convergence}.
\subsection{Related works}

Recently there are a lot of papers discussing various
atomistic/continuum coupling strategies as summarized in the recent
reviews~\cites {RuddBroughton:2000, CurtinMiller:2003,
  MillerTadmor:2009, HMMreview}, we will only mention some of the
works that are closely related to ours and refer the readers to these
reviews and the references therein.

The hybrid method resembles several methods in the literature. The
most closely related method is the quasicontinuum (QC)
method~\cites{TadmorOrtizPhillips:1996, KnapOrtiz:2001}, which is
among the most popular methods for modeling the mechanical
deformation of crystalline solids.  The QC method contains following
ingredients: decomposition of the whole domain into atomistic
and continuum regions, with the defects covered by the atomistic
regions; degree reduction by adaptive selection of representative
atoms (rep-atoms), with fewer atoms selected in regions with smooth
deformation; and the application of the Cauchy-Born approximation in
the continuum region to reduce the complexity involved in computing
the total energy of the system.

Both the proposed method and QC method couple atomistic models with
nonlinear Cauchy-Born elasticity model. In some sense, the proposed
method can be viewed as a smoothened modification of the force-based
QC method. Indeed, the original force-based QC method amounts to
take $\varrho$ to be a characteristic function (so that there is no
buffer region). The force-based QC is free of ghost force, and it
was proven in~\cites{Ming:2008, DobsonLuskinOrtner:2010a} that, for
one-dimensional problem, the force-based method converges
quadratically. However, its convergence behavior remains open for
high dimensional problem. As will be proved later in the paper, the
proposed method is stable and also converges quadratically in three
dimension. For the understanding of the original force-based QC,
this work may also provide some new tools and insights.

The Arlequin method~\cites{BenDhia:1998, Bauman:2008} and the bridging
domain method~\cite{BelytschkoXiao:2003} also adopt a smooth
transition between atomistic and continuum regions. The difference
with the proposed scheme is however these methods are energy-based, so
that the mixing is done at the energy level, while the current method
is force-based. Moreover, these two methods enforce the consistency
between the atomistic and continuum regions by imposing certain
constraints, while there is no such constraints in our method.  These
methods are suffered from ghost force problems as shown
in~\cite{MillerTadmor:2009}, while the proposed method is consistent
at the interface.

The proposed method also shares certain common traits with the
concurrent AtC coupling method (AtC) proposed
in~\cite{BadiaBochvLehoucqParksFishNuggehallyGunzburger:2007}. The AtC
method also uses a smooth transition between atomistic and continuum
regions and is force-based. However, the proposed method differs from
AtC in the following aspects:
(1) our method employs Cauchy-Born elasticity while AtC uses linear
elasticity and (2) our method is free of ghost force while AtC is
plagued by ghost force as demonstrated in~\cite{MillerTadmor:2009}.

Most of the analysis of these multiscale methods limits to the
quasicontinuum method.  In~\cite{EMing:2007}, the Cauchy-Born rule
for crystalline solids is verified under sharp stability conditions.
In the language of QC, the authors in~\cite{EMing:2007} actually
proved the convergence of local QC (the whole computational domain
is treated as local region). Explicit convergence rate for the local
QC can be found in~\cites{EMing:2004,EMing:2005}.

For the QC method couples together atomistic and continuum models
(nonlocal QC method in short), the error estimate can be found
in~\cites{MingYang:2009, DobsonLuskin:2009b} and the references
therein. All these works dealt only with the one dimensional problem,
and moreover, except \cite{MingYang:2009}, the analysis was limited to
quadratic potential models, so that the system is linear.

To the best of the authors' knowledge, there is no analysis for the
nonlocal QC method or other coupling schemes for high-dimensional
problems with general potential (usually, many-body potential
function). The main difficulties lie in the analysis of the
consistency and stability. For one-dimensional problem, the lattice
structure is very simple and the pairwise potential function can be
handled by a direct calculation. However, such an approach cannot be
easily extended to high-dimensional problem with general potential
because the lattice structure and the potential function for
high-dimensional problem is much more involved.  One of the main
contributions of the current paper is the development of general tools
for the analysis of consistency and stability.

Finally, we remark that in this work the analysis of the proposed
method, especially the stability analysis, is based on analysis of
finite difference schemes. The readers might wonder why the analysis
is \emph{not} done in the framework of finite element method, as
after all, we are dealing with static problems, the systems to be
solved are ``elliptic''; and moreover, the continuum region is
discretized by finite element method. The reason actually lies in
the atomistic part, since the force balance equations derived from
energy of discrete lattice systems are intrinsically of finite
difference type. To the best of our knowledge, there has not been
yet a successful way to put the atomistic equations into the
framework of finite element analysis. Therefore, to be consistent,
we view the finite element approximation in the continuum region
also as a finite difference approximation. The proof hence relies on
the analysis of finite difference schemes. This may give a
reminiscence of the early history about finite element analysis,
during when the finite element method was also analyzed in the
framework of finite difference schemes~\cite{StrangFix:1973}. Since
the theory of adaptive mesh is well-established for finite
element method, it is an interesting question whether one can
adopt the finite element analysis framework to analyze these
multiscale coupling methods.
\section{Consistency}\label{sec:consistency}

We study the consistency of the force-based hybrid method in this
section. The key is the following lemma, which is a refined version
of~\cite{EMing:2007}*{Lemma 5.1}.
\begin{lemma}[Consistency of Cauchy-Born rule]\label{lem:consCB}
  For any $y=x+u(x)$ with $u$ smooth, we have
  \begin{equation}\label{cons:eq1}
    \norm{\mc{F}_{\at}[y] - \mc{F}_{\CB}[y]}_{L^{\infty}_{\veps}} \leq
    C \veps^2 \norm{u}_{W^{16, \infty}},
  \end{equation}
  where the constant $C$ depends on $V$ and $\norm{u}_{L^\infty}$, but
  is independent of $\veps$.
\end{lemma}

\begin{remark}
  The consistency estimate is presented in the form of
  \eqref{cons:eq1} for later use in the proof of
  Proposition~\ref{prop:Hcons}. A bound involves less order of
  derivatives of $u$ is possible, in fact, it is not
  hard to see from the proof that we have
  \begin{equation}\label{cons:lowdiff}
    \norm{\mc{F}_{\at}[y] - \mc{F}_{\CB}[y]}_{L^{\infty}_{\veps}} \leq
    C  \veps^2,
  \end{equation}
  where $C$ depends on $V$ and $\norm{u}_{W^{6, \infty}}$. The price
  is however the dependence of $C$ on $\norm{u}_{W^{6, \infty}}$ is
  nonlinear.
\end{remark}

\begin {proof}
For any $x\in\Omega_{\veps}$, and for $i=1,2$, Taylor expansion at
$x$ gives
\[
\ffd{s_i}y(x) = \na_{s_i}^1[y](x)+\veps\na_{s_i}^2[y](x)+\veps^2
R_{2,s_i}[y](x),
\]
where, for convenience, we have introduced the short-hands for the
Taylor series and its remainder:
\begin {align*}
  & \na_{s_i}^j[y](x) = \dfrac{1}{j!}(s_i\cdot\na)^j y(x), \\
  & R_{k,s_i}[y](x) =\int_0^1(k+1)(1-t)^k\na_{s_i}^{k+1} y(x+\veps
  ts_i)\dt,\quad k\in\NN,
\end {align*}
provided that the terms on the right hand side are well defined.
Obviously, we may write
\begin {equation}\label {operexpan}
\ffd{s_i}=\na_{s_i}^1+\veps\na_{s_i}^2+\veps^2 R_{2,s_i},\quad
\bfd{s_i}=\na_{s_i}^1-\veps\na_{s_i}^2-\veps^2 R_{2,-s_i}.
\end {equation}

For $i=1,2,3$ and $t \in [0,1]$, let
\begin{align*}
  F_i(t) =\pa_i V_{(s_1,s_2)}\Bigl(
  & \abs{t\ffd{s_1}y(x)+(1-t)(s_1\cdot\na)y(x)}^2,\\
  & \quad \abs{t\ffd{s_2}y(x)
    +(1-t)(s_2\cdot\na)y(x)}^2,\\
  & \quad \bigl\langle t\ffd{s_1}y(x)+(1-t)(s_1\cdot\na)y(x), \\
  & \qquad\qquad t\ffd{s_2}y(x) +(1-t)(s_2\cdot\na)y(x)\bigr\rangle
  \Bigr).
\end{align*}
Using Taylor expansion, we get
\begin{equation}\label{eq:expandF}
  F_i(1)=F_i(0)+F^{\prime}_i(0)+R_1[F_i](0).
\end{equation}
Here for $F_i: [0, 1] \to \RR$, we have introduced a similar
short-hand for the remainder
\[
R_k[F_i](0)=\int_0^1 \frac{(1-t)^k}{k!}\na^{k+1}F_i(t) \ud t.
\]

Notice that by definition we have
\begin{align*}
  F_i(1) & = \partial_i V_{(s_1, s_2)}\bigl(\abs{\ffd{s_1}y(x)}^2,
  \abs{\ffd{s_2}y(x)}^2, \inner{\ffd{s_1}y(x)}{\ffd{s_2}y(x)} \bigr) \\
  & = \partial_i V_{(s_1, s_2)}[y](x); \\
  F_i(0) & = \partial_i V_{(s_1, s_2)}\bigl(\abs{(s_1\cdot\na)y(x)}^2,
  \abs{(s_2\cdot\na)y(x)}^2,
  \inner{(s_1\cdot\na)y(x)}{(s_2\cdot\na)y(x)} \bigr) \\
  & = \partial_i W_{(s_1, s_2)}(\nabla u(x)).
\end{align*}
Therefore, we can rewrite \eqref{eq:expandF} as
\begin{equation}\label{force-expan}
  \begin{aligned}
    \pa_iV_{(s_1,s_2)}[y](x)&=\pa_iW_{(s_1,s_2)}(\na u(x))+\veps
    a_j\pa_{ij}W_{(s_1,s_2)}(\na u(x))\\
    &\quad+\Lr{\veps^2
      b_j\pa_{ij}W_{(s_1,s_2)}(\na u(x))+R_1[F_i](0)}\\
    &\quad+\veps^3c_j\pa_{ij}W_{(s_1,s_2)}(\na u(x))\\
    &\equiv\mc{Q}_{i, (s_1, s_2)}[\na u](x),
  \end{aligned}
\end{equation}
where for $j = 1, 2, 3$,
\begin {align*}
a_j&=2\inner{(s_j\cdot\na)y}{\na_{s_j}^2[y]}(1-\delta_{j3})\\
&\quad+\Lr{\inner{(s_1\cdot\na)y}{\na_{s_2}^2[y]}
+\inner{(s_2\cdot\na)y}{\na_{s_1}^2[y]}}\delta_{j3},\\
b_j&=2\inner{(s_j\cdot\na)y}{\na_{s_j}^3[y]}(1-\delta_{j3})\\
&\quad+\Lr{\inner{(s_1\cdot\na)y}{\na_{s_2}^3[y]}
+\inner{(s_2\cdot\na)y}{\na_{s_1}^3[y]}}\delta_{j3},\\
c_j&=2\inner{(s_j\cdot\na)y}{R_{2,s_j}[y]}(1-\delta_{j3})\\
&\quad+\Lr{\inner{(s_1\cdot\na)y}{R_{2,s_2}[y]}
+\inner{(s_2\cdot\na)y}{R_{2,s_1}[y]}}\delta_{j3}.
\end {align*}

Substituting the equations~\eqref {operexpan} into
$\mc{F}_{\at}[y](x)$, we obtain
\begin {align*}
  \mc{F}_{\at}[y]=& \\
  \sum_{(s_1, s_2) \in S}
  &(\na_{s_1}^1-\veps\na_{s_1}^2-\veps^2R_{2,-s_1})
  \Bigl\{2\pa_1V_{(s_1,s_2)}[y](\na_{s_1}^1+\veps\na_{s_1}^2+\veps^2R_{2,s_1})[y]\\
  &\phantom{(\na_{s_1}^1-\veps\na_{s_1}^2+\veps^2R_{2,-s_1})}\qquad
  +\pa_3V_{(s_1, s_2)}[y](\na_{s_2}^1+\veps\na_{s_2}^2+\veps^2R_{2,s_2})[y]\Bigr\}\\
  +&(\na_{s_2}^1-\veps\na_{s_2}^2-\veps^2R_{2,-s_2})\Bigl\{2\pa_2
  V_{(s_1,s_2)}[y](\na_{s_2}^1+\veps\na_{s_2}^2+\veps^2R_{2,s_2})[y]\\
  &\phantom{(\na_{s_2}^1-\veps\na_{s_2}^2+\veps^2R_{2,-s_2})}\qquad
  +\pa_3V_{(s_1,s_2)}[y](\na_{s_2}^1+\veps\na_{s_2}^2+\veps^2R_{2,s_2})[y]\Bigr\}.
\end {align*}
Next substituting~\eqref {force-expan} into the above equation, we have
\begin {align*}
  \mc{F}_{\at}[y](x)= &\\
  \sum_{(s_1, s_2) \in S} &
  (\na_{s_1}^1-\veps\na_{s_1}^2-\veps^2R_{2,-s_1}) \Bigl\{2\mc{Q}_{1,
    (s_1, s_2)}[\na u](\na_{s_1}^1+\veps\na_{s_1}^2+\veps^2
  R_{2,s_1})[y]\\
  &\phantom{(\na_{s_1}^1-\veps\na_{s_1}^2+\veps^2R_{2,-s_1})}\
  +\mc{Q}_{3, (s_1, s_2)}[\na u]
  (\na_{s_2}^1+\veps\na_{s_2}^2+\veps^2 R_{2,s_2})[y]\Bigr\}\\
  + & (\na_{s_2}^1-\veps\na_{s_2}^2-\veps^2R_{2,-s_2})
  \Bigl\{2\mc{Q}_{2, (s_1, s_2)}[\na u](\na_{s_2}^1+\veps\na_{s_2}^2
  +\veps^2 R_{2,s_2})[y]\\
  &\phantom{(\na_{s_2}^1-\veps\na_{s_2}^2+\veps^2R_{2,-s_2})}\
  +\mc{Q}_{3, (s_1, s_2)}[\na u](\na_{s_1}^1+\veps\na_{s_1}^2+\veps^2
  R_{2,s_1})[y]\Bigr\}.
\end {align*}
Collecting the terms of the same order, we get
\begin {equation}\label {opexpansion:eq}
\mc{F}_{\at}[y](x)= \mc{L}_0[u](x)
+\veps\mc{L}_1[u](x)+\veps^2\mc{L}_2[u](x)+\mc{O}(\veps^3).
\end {equation}
If we change $\veps$ to $-\veps$, the left-hand side of~\eqref
{opexpansion:eq} is invariant, then the terms of odd power of
$\veps$ in the right-hand side of~\eqref {opexpansion:eq}
automatically vanishes. Therefore, we have
\[
\mc{F}_{\at}[y](x)=\mc{L}_0[u](x)
+\veps^2\mc{L}_2[u](x)+\mc{O}(\veps^4).
\]

The explicit form of $\mc{L}_0$  can be written as
\begin {align*}
  \mc{L}_0[u](x) &=
  -2(s_1\cdot\na)\bigl[(s_1+(s_1\cdot\na)u)\pa_1W_{(s_1,s_2)}(\na u(x))\bigr]\\
  &\quad-(s_1\cdot\na)\bigl[(s_2+(s_2\cdot\na)u)
  \pa_3W_{(s_1,s_2)}(\na u(x))\bigr]\\
  &\quad-2(s_2\cdot\na)\bigl[(s_2+(s_2\cdot\na)u)
  \pa_2W_{(s_1,s_2)}(\na u(x))\bigr]\\
  &\quad-(s_2\cdot\na)\bigl[(s_1+(s_1\cdot\na)u)\pa_3W_{(s_1,s_2)}(\na
  u(x))\bigr].
\end {align*}
We see that $\mc{L}_0$ is the same as the operator that appears in
the Euler-Lagrangian equation of~\eqref {cb:min}.

The proof of that $\mc{L}_2$ is of divergence form is similar.
Actually, $\mc{L}_2$ is a quasilinear operator, which actually
counts for the linear dependence on $\norm{u}_{W^{16, \infty}}$ on
the right-hand side of~\eqref{cons:eq1}. To prove~\eqref{cons:eq1},
it remains to estimate terms of $\Or(\veps^2)$, which is a
combination of terms of the form: for $\al,\beta=1,2$,
\begin {align*}
&\na_{s_\al}^k\Lr{\pa_iW_{(s_1,s_2)}(\na u)\na_{s_\beta}^l
u},\quad l+k=4,l,k\in\NN,\\
&\na_{s_\al}^k\Lr{a_j\pa_{ij}W_{(s_1,s_2)}(\na u)\na_{s_\beta}^l
u},\quad l+k=3,l,k\in\NN,\\
&\na_{s_\al}^1\Lr{b_j\pa_{ij}W_{(s_1,s_2)}(\na u)\na_{s_\beta}^1
u+R_1[F_i](0)\na_{s_{\beta}}^1u}.
\end {align*}
We only give the estimate for the first term, and the other two can be
bounded similarly. Due to chain rule and to Leibniz's rule,
$\na_{s_\al}^k\Lr{\pa_iW_{(s_1,s_2)}(\na u)\na_{s_\beta}^l u}$ is a
linear combination of the form
\begin {align*}
T&=\prod_{i=1}^3\Lr{\dfrac{\pa}{\pa x_i}}^{\sgn{\delta_i}}\pa_i
W_{(s_1,s_2)}(\na
u)\\
&\phantom{=\prod_{i=1}^3}\quad\x
(s_\al\cdot\na)^{\ga_1}P_{\delta_1}(s_\al\cdot\na)^{\ga_2}P_{\delta_2}
(s_\al\cdot\na)^{\ga_3}P_{\delta_3}(s_{\beta}\cdot\na)^{4-\abs{\ga}}u,
\end {align*}
where $\ga\in\NN^3$ are multiindecies with
$\abs{\ga}=\sum_{i=1}^3\abs{\ga_i}$ and $\abs{\ga}\le 3$. Here
\[
P_1=\abs{s_1+(s_1\cdot\na)u}^2,\quad
P_2=\abs{s_2+(s_2\cdot\na)u}^2,\quad
P_3=\inner{s_1+(s_1\cdot\na)u}{s_2+(s_2\cdot\na)u}.
\]
Using chain rule once again, we get, for $i=1,2,3$,
\begin {align*}
\norm{(s_\al\cdot\na)P_i}_{L^\infty}&\le
C(s_\al)(1+\norm{\na u}_{L^\infty})\norm{\na^2 u}_{L^\infty},\\
\norm{(s_\al\cdot\na)^2P_i}_{L^\infty}&\le
C(s_\al)\Lr{(1+\norm{\na u}_{L^\infty})\norm{\na^3 u}_{L^\infty}
+\norm{\na^2 u}_{L^\infty}^2},\\
\norm{(s_\al\cdot\na)^3P_i}_{L^\infty}&\le C(s_\al)\Lr{(1+\norm{\na
u}_{L^\infty})\norm{\na^4 u}_{L^\infty}+\norm{\na^2
u}_{L^\infty}^2\norm{\na^3 u}_{L^\infty}^2}.
\end {align*}
Using Gagliardo-Nirenberg inequality~\cite{Nirenberg: 1959},
\[
\norm{\na^j u}_{L^\infty}\le C\norm{\na^m
u}_{L^\infty}^{j/m}\norm{u}_{L^\infty}^{1-j/m},\quad 0<j<m,
\]
we have
\[
\norm{(s_\al\cdot\na)^kP_i}_{L^\infty}\le C(s_\al)\Lr{\norm{
u}_{L^\infty}\norm{\na^{k+2} u}_{L^\infty}+\norm{\na^{k+1}
u}_{L^\infty}}.
\]
Using the above inequality, we conclude
\begin {align*}
\norm{T}_{L^\infty}&\le C\max_{2\le\abs{\ga}\le
4}\norm{\pa_{\ga}W_{(s_1,s_2)}(\na
u)}_{L^\infty}\norm{\na^{4-\abs{\ga}} u}_{L^\infty}\\
&\quad\x\Biggl\{
(1+\norm{u}_{L^\infty}^3)\prod_{i=1}^3\norm{\na^{{\ga_i}+2}u}_{L^\infty}
+\prod_{i=1}^3\norm{\na^{{\ga_i}+1}u}_{L^\infty}\\
&\phantom{\x\Biggl\{}\qquad+(1+\norm{u}_{L^\infty}^2)\sum_{i,j,k=1}^3\norm{\na^{{\ga_i}+2}u}_{L^\infty}
\norm{\na^{{\ga_j}+2}u}_{L^\infty}\norm{\na^{{\ga_k}+1}u}_{L^\infty}\\
&\phantom{\x\Biggl\{}\qquad+(1+\norm{u}_{L^\infty})\sum_{i,j,k=1}^3\norm{\na^{{\ga_i}+2}u}_{L^\infty}
\norm{\na^{{\ga_j}+1}u}_{L^\infty}\norm{\na^{{\ga_k}+1}u}_{L^\infty}\Biggr\}.
\end {align*}
Invoking Gagliardo-Nirenberg inequality again, we obtain
\begin {align*}
\norm{T}_{L^\infty}&\le
C(\norm{u}_{L^\infty}^3+\norm{u}_{L^\infty}^6)\norm{\na^{10}u}_{L^\infty}\\
&\quad+C(\norm{u}_{L^\infty}^3+\norm{u}_{L^\infty}^5)\norm{\na^9u}_{L^\infty}\\
&\quad+C(\norm{u}_{L^\infty}^3+\norm{u}_{L^\infty}^4)\norm{\na^8u}_{L^\infty}\\
&\quad+C\norm{u}_{L^\infty}^3\norm{\na^7u}_{L^\infty}\\
&\le C\sum_{i=3}^6\norm{u}_{L^\infty}^i\norm{u}_{W^{10,\infty}}.
\end {align*}
Proceeding along the same line, we can obtain the similar bounds for
the higher-order terms, while $\norm{u}_{W^{16,\infty}}$ arises from
the following term
\[
R_{\al,-s_\al}\Lr{\pa_i W_{(s_1,s_2)}(\na u)R_{\beta,s_{\beta}}[y]}.
\]
Summing up all terms of $\Or{(\veps^2)}$, we get~\eqref{cons:eq1}.

\end{proof}
\begin{coro}[Consistency of finite element discretization]\label{coro:consFE}
  For any $y=x+u(x)$ with $u$ smooth, we have
  \begin{equation*}
    \norm{\mc{F}_{\veps}[y] - \mc{F}_{\CB}[y]}_{L^{\infty}_{\veps}}
    \leq C\veps^2 \norm{u}_{W^{16, \infty}},
  \end{equation*}
  where the constant $C$ depends on $V$ and $\norm{u}_{L^\infty}$, but
  is independent of $\veps$.
\end{coro}

\begin{proof}
  The corollary follows Lemma~\ref{lem:consCB} by the observation that
  we can view the energy functional of the finite element
  discretization as a particular choice of atomistic potential energy.

  To be more concrete, let us consider the case $d = 2$, so that each
  element $T \in \mc{T}_{\veps}$ has three vertices. It is
  straightforward to extend the argument below to higher dimensions,
  with certain complication of notations.

  Let $y_{\veps} \in \wt{X}_{\veps}$ be the approximation of $y$ so
  that $y_{\veps}(x) = y(x)$ for any $x \in \Omega_{\veps}$. Let
  $u_{\veps} = y_{\veps} - x$. Obviously, we have $u_{\veps}(x) =
  u(x)$ for any $x \in \Omega_{\veps}$.

  Now, for each $T \in \mc{T}_{\veps}$, $\nabla u_{\veps}
  \vert_{T}$ is a linear function of $y_{\veps}$ on the
  vertices of $T$. Denote the three vertices of $T$ as
  $x_0, x_1, x_2$, and $s_1 = (x_1 - x_0)/\veps$, $s_2 = (x_2 - x_0) /
  \veps$, then $\nabla u_{\veps}\vert_{T}$ is the solution of
  the linear system
  \begin{equation*}
    \begin{cases}
      s_1 + s_1 A = D_{\veps, s_1}^+ y_{\veps}(x_0), \\
      s_2 + s_2 A = D_{\veps, s_2}^+ y_{\veps}(x_0).
    \end{cases}
  \end{equation*}
  Therefore, let us denote
  \begin{equation*}
    \nabla u_{\veps} \vert_{T} = A_{(s_1, s_2)}(y_{\veps}(x_0)/\veps,
    y_{\veps}(x_1)/\veps, y_{\veps}(x_2) / \veps)
  \end{equation*}
  as the solution of the above system. Notice that due to linearity,
  the map $A_{(s_1, s_2)}$ is independent of $\veps$.  Hence, for $x
  \in T$, we can write
  \begin{equation}\label{eq:Wcbpotential}
    \begin{aligned}
      W_{\CB}(\nabla u_{\veps}(x)) & = W_{\CB}\bigl(A_{(s_1,
        s_2)}(y_{\veps}(x_0)/\veps,
      y_{\veps}(x_1)/\veps, y_{\veps}(x_2) / \veps)\bigr) \\
      & = W_{\FE, (s_1, s_2)}(y_{\veps}(x_0)/\veps,
      y_{\veps}(x_1)/\veps, y_{\veps}(x_2) / \veps),
    \end{aligned}
  \end{equation}
  where $W_{\FE, (s_1, s_2)} \equiv W_{\CB}\circ A_{(s_1, s_2)}$. Denote
  $S_{\FE}$ as the set of all pairs $(s_1, s_2)$ such that $\{x_0,
  x_0+\veps s_1, x_0 + \veps s_2 \}$ forms the vertices of an element
  $T \in \mc{T}_{\veps}$ containing $x_0$ (it is easy to see that
  $S_{\FE}$ is independent of $\veps$). Then, using
  \eqref{eq:Wcbpotential}, we have
  \begin{equation*}
    \begin{aligned}
      \int_{\Omega} & W_{\CB}(\nabla u_{\veps}(x)) \\
      & = \sum_{T
        \in \mc{T}_{\veps}} \abs{T} W_{\CB}(\nabla u_{\veps} \vert_T) \\
      & = \frac{1}{3!}\sum_{x \in \Omega_{\veps}} \sum_{(s_1, s_2) \in
        S_{\FE}} \veps^d \abs{T_{(s_1, s_2)}} W_{\FE, (s_1,
        s_2)}\left(\frac{y_{\veps}(x)}{\veps}, \frac{y_{\veps}(x +
          \veps s_1)}{\veps}, \frac{y_{\veps}(x + \veps
          s_2)}{\veps}\right) \\
      & = \frac{1}{3!}\veps^d \sum_{x \in \Omega_{\veps}} \sum_{(s_1,
        s_2) \in S_{\FE}} V_{\FE, (s_1,
        s_2)}\left(\frac{y_{\veps}(x)}{\veps}, \frac{y_{\veps}(x +
          \veps s_1)}{\veps}, \frac{y_{\veps}(x + \veps
          s_2)}{\veps}\right),
    \end{aligned}
  \end{equation*}
  where $V_{\FE, (s_1, s_2)} =\abs{T_{(s_1, s_2)}} W_{\FE, (s_1,
    s_2)}$ and $T_{(s_1, s_2)}$ is the triangle formed by vectors
  $s_1$ and $s_2$. This indicates that we can view the energy
  functional in the finite element discretization as a particular
  atomistic potential model, given by three body interactions $V_{\FE,
    (s_1, s_2)}$, by identifying the value of $y$ on nodes as the
  deformed atom positions.

  It is immediately clear that the Cauchy-Born energy density
  corresponding to the atomic potential constructed is just
  $W_{\CB}$. Indeed, for a homogenously deformed system with
  deformation gradient $A$, by definition, the energy of the system is
  just $W_{\CB}(A) \abs{\Omega}$, and hence the Cauchy-Born energy
  density is given again by $W_{\CB}(A)$.

  With this viewpoint of the finite element discretization as an
  atomic potential, we obtain the conclusion as an immediate corollary
  of Lemma~\ref{lem:consCB}.
\end{proof}

\begin{coro}[Local truncation error]\label{coro:consqc}
  For any $y=x+u(x)$ with $u$ smooth, we have
  \begin{equation}\label{eq:compareqcat}
    \norm{\mc{F}_{\qc}[y] - \mc{F}_{\at}[y]}_{L^{\infty}_{\veps}} \leq C \veps^2
    \norm{u}_{W^{16, \infty}},
  \end{equation}
  and
  \begin{equation}\label{eq:compareqcCB}
    \norm{\mc{F}_{\qc}[y] - \mc{F}_{\CB}[y]}_{L^{\infty}_{\veps}} \leq C \veps^2
    \norm{u}_{W^{16, \infty}},
  \end{equation}
  where the constant $C$ depends on $V$ and $\norm{u}_{L^\infty}$, but
  is independent of $\veps$.
\end{coro}

\begin{proof}
  The inequality \eqref{eq:compareqcat} follows from
  Lemma~\ref{lem:consCB}, Corollary~\ref{coro:consFE}, and
  \[
    \begin{aligned}
      \norm{\mc{F}_{\qc}[y] - \mc{F}_{\at}[y]}_{L^{\infty}_{\veps}} &
      = \norm{ \rho(x) (\mc{F}_{\veps}[y](x) - \mc{F}_{\at}[y](x))
      }_{L^{\infty}_{\veps}} \\
      & \leq \norm{ \mc{F}_{\veps}[y] - \mc{F}_{\at}[y]}_{L^{\infty}_{\veps}} \\
      & \leq \norm{ \mc{F}_{\at}[y] -
        \mc{F}_{\CB}[y]}_{L^{\infty}_{\veps}} + \norm{
        \mc{F}_{\veps}[y] - \mc{F}_{\CB}[y]}_{L^{\infty}_{\veps}},
    \end{aligned}
  \]
  where we have used $\varrho(x) \in [0, 1]$.  Similarly,
  \eqref{eq:compareqcCB} follows from Lemma~\ref{lem:consCB},
  Corollary~\ref{coro:consFE}, and
  \[
    \begin{aligned}
      \norm{\mc{F}_{\qc}[y] - \mc{F}_{\CB}[y]}_{L^{\infty}_{\veps}} &
      \leq \norm{ \rho(x) (\mc{F}_{\veps}[y](x) - \mc{F}_{\CB}[y](x))
      }_{L^{\infty}_{\veps}} \\
      & \qquad + \norm{(1 - \rho(x)) (\mc{F}_{\at}[y](x)
        - \mc{F}_{\CB}[y](x))}_{L^{\infty}_{\veps}}  \\
      & \leq \norm{ \mc{F}_{\veps}[y] - \mc{F}_{\CB}[y]
      }_{L^{\infty}_{\veps}} + \norm{ \mc{F}_{\at}[y] -
        \mc{F}_{\CB}[y]}_{L^{\infty}_{\veps}}.
    \end{aligned}
  \]

\end{proof}
\section{Regularity estimate}\label{sec:regularity}

To analyze the stability property of the proposed force-based hybrid
method, we use the framework of pseudo-difference operators
\cites{Thomee:1964, LaxNirenberg:1966}. In this section, we will
establish regularity estimate Theorem~\ref{thm:regularity} for the
force-based hybrid method. This will be one of the key ingredients
used to prove stability estimate in the next section.

We study the linearized operator of $\mc{F}_{\qc}$.  Let us denote
$\mc{H}_{\qc}[u]$ the linearization of $\mc{F}_{\qc}$ at state $u$:
\[
  \mc{H}_{\qc}[u] = \frac{\delta \mc{F}_{\qc}}{\delta y}
  \bigg\vert_{y = x + u},
\]
so that $\mc{H}_{\qc}[u]$ is a linear operator on lattice functions
$w$, given by
\[
  \mc{H}_{\qc}[u] w = \lim_{t\to0} \frac{\partial \mc{F}_{\qc}[x + u + t w]}
  {\partial t}.
\]
It is convenient to rewrite $\mc{H}_{\qc}$ in the form of a
pseudo-difference operator as
\[
  \mc{H}_{\qc}[u] = \sum_{\mu \in \mc{A}} h_{\qc}[u](x, \mu) T^{\mu},
\]
where the coefficient $h_{\qc}[u](x,\mu)$ is a $d\times d$ (probably
asymmetric) matrix for each $x$ and $\mu \in \mc{A}$, given by
\begin{equation}\label{eq:defhqc}
  (h_{\qc}[u])_{\alpha\beta}(x, \mu) = \frac{\partial
    ( \mc{F}_{\qc}[y] )_{\alpha}(x)}
  {\partial (T^{\mu} y)_{\beta}(x)} \bigg\vert_{y = x + u},
\end{equation}
where $\alpha, \beta = 1, \cdots, d$ are indices. Here $\mc{A}$ is
range of the pseudo-difference stencil (note that $0 \in \mc{A}$),
which is finite by assumptions. By the definition of $\mc{F}_{\qc}$,
we have
\begin{equation}\label{eq:linearh}
  h_{\qc}[u](x, \mu) = (1 - \varrho(x)) h_{\at}[u](x, \mu)
  + \varrho(x) h_{\veps}[u](x, \mu),
\end{equation}
where $h_{\at}[u]$ and $h_{\veps}[u]$ are given by similar equations
as \eqref{eq:defhqc} by replacing $\mc{F}_{\qc}$ to $\mc{F}_{\at}$ and
$\mc{F}_{\veps}$ respectively.

Define $\wt{h}_{\qc}[u](x, \xi)$ as the symbol of the
pseudo-difference operator $\mc{H}_{\qc}[u]$ given as
\[
  \wt{h}_{\qc}[u](x, \xi) = \sum_{\mu\in\mc{A}} h_{\qc}[u](x, \mu)
  \exp(\I \veps \sum_j \mu_j a_j \cdot \xi) \qquad \text{for } \xi \in
  \LL^{\ast}_{\veps},
\]
and similarly for $\wt{h}_{\veps}[u]$ and $\wt{h}_{\at}[u]$.  By
definition, we have for any $x \in \Omega_{\veps}$,
\[
  (\mc{H}_{\qc}[u] e_k e^{\I x\cdot \xi})_j(x) =
  (\wt{h}_{\qc}[u])_{jk}(x, \xi) e^{\I x\cdot \xi},
\]
for $1 \leq j,k \leq d$ and similarly for $\wt{h}_{\veps}[u]$ and
$\wt{h}_{\at}[u]$. Here $\{e_k\}$ are the canonical basis of
$\RR^d$. It is also clear that~\eqref{eq:linearh} implies
\begin{equation}\label{eq:linearwth}
  \wt{h}_{\qc}[u](x, \xi) = (1 - \varrho(x)) \wt{h}_{\at}[u](x, \xi)
  + \varrho(x) \wt{h}_{\veps}[u](x, \xi).
\end{equation}

In the case that we linearize around the equilibrium state $u = 0$,
we will simplify the notation as
\begin{equation*}
  \mc{H}_{\qc} = \mc{H}_{\qc}[0], \quad
  h_{\qc} = h_{\qc}[0], \quad
  \wt{h}_{\qc} = \wt{h}_{\qc}[0],
\end{equation*}
and similarly for those defined for atomistic model and finite element
discretization. We observe that by the translation invariance of the
total energy $I_\at$ at the state $u = 0$,
\[
  h_{\at}(x, \mu) = h_{\at}(\mu), \quad h_{\veps}(x, \mu)
  = h_{\veps}(\mu).
\]
The coefficients are independent of position $x$, and hence similarly
for $\wt{h}_{\at}$ and $\wt{h}_{\veps}$.

We also denote $\mc{H}_{\CB}$ as the linearization of $\mc{F}_{\CB}$
at the equilibrium state $u = 0$, and define $\wt{h}_{\CB} =
\wt{h}_{\CB}(x, \xi)$ as its symbol. Note that due to the periodic
boundary condition assumed on $\Omega$, $\xi$ here only takes value
in $\LL^{\ast}$. Again, due to the translation invariance of the
total energy, $\wt{h}_{\CB}$ is independent of $x$.

Let us start the analysis with the operator $\mc{H}_{\qc}$. First,
we show that the matrix $\wt{h}_{\qc}$ is Hermitian.
\begin{lemma}
  The matrices $\wt{h}_{\at}(\xi)$, $\wt{h}_{\veps}(\xi)$ and hence
  $\wt{h}_{\qc}(x, \xi)$ are Hermitian for any $\veps > 0$, $x
    \in \Omega_{\veps}$ and $\xi \in \LL_{\veps}^{\ast}$.
\end{lemma}
\begin{proof}
  It suffices to prove the result for $\wt{h}_{\at}(\xi)$, as the
  argument for $\wt{h}_{\veps}(\xi)$ is the same and the
  conclusion for $\wt{h}_{\qc}(x, \xi)$ follows immediately from
  \eqref{eq:linearwth}.

  Since $(\mc{F}_{\at}[y])_{\alpha}(x) = - \partial I_{\at}[y]
  / \partial y_{\alpha}(x)$, we have
  \begin{align*}
    (h_{\at})_{\alpha\beta}(\mu) & = - \frac{\partial^2 I_{\at}[y]}
    {\partial y_{\alpha}(x) \partial (T^{\mu} y)_{\beta}(x)} \bigg\vert_{y = x} \\
    & = - \frac{\partial^2 I_{\at}[y]} {\partial
      y_{\alpha}(x) \partial y_{\beta}(x + \veps\mu_j a_j)}
    \bigg\vert_{y = x} \\
    & = - \frac{\partial^2 I_{\at}[y]} {\partial
      (T^{-\mu}y)_{\alpha}(x+\veps \mu_j a_j) \partial y_{\beta}(x +
      \veps\mu_j a_j)} \bigg\vert_{y = x} \\
    & = - \frac{\partial^2 I_{\at}[y]} {\partial
      (T^{-\mu}y)_{\alpha}(x) \partial y_{\beta}(x)} \bigg\vert_{y =
      x} = (h_{\at})_{\beta\alpha}(-\mu),
  \end{align*}
  where the last line follows from translational invariance of the
  unperturbed system. Therefore,
  \begin{align*}
    (\wt{h}_{\at})_{\alpha\beta}(\xi) & = \sum_{\mu}
    (h_{\at})_{\alpha\beta}(\mu) \exp(\I \veps \sum_j \mu_j a_j \cdot \xi) \\
    & = \sum_{\mu} (h_{\at})_{\beta\alpha}(-\mu)
    \exp(\I \veps \sum_j (-\mu_j) a_j \cdot (-\xi)) \\
    & = \biggl( \sum_{\mu} (h_{\at})_{\beta\alpha}(-\mu) \exp(\I \veps
    \sum_j (-\mu_j) a_j \cdot \xi) \biggr)^{\ast} =
    (\wt{h}_{\at})_{\beta\alpha}^{\ast}(\xi),
  \end{align*}
  for any $\xi \in \LL_{\veps}^{\ast}$, where we have used the fact
  that $h_{\at}$ are real matrices. This proves the Lemma.
\end{proof}

We make the following stability assumptions about the atomistic
potentials, the finite element discretization of the
Cauchy-Born elasticity model:
\begin{assump}\label{assump:stabatom}
  $\wt{h}_{\at}(\xi)$ is positive definite and there exists
  $a_{\at}>0$ such that for any $\veps>0$ and any $\xi \in
  \LL_{\veps}^{\ast}$,
  \begin{equation*}
    \det \wt{h}_{\at}(\xi) \geq a_{\at} \Lambda_{0, \veps}^{2d}(\xi).
  \end{equation*}
\end{assump}
\begin{assump}\label{assump:stabCB}
  $\wt{h}_{\veps}(\xi)$ is positive definite and there exists $a>0$
  such that for any $\veps>0$ and any $\xi \in \LL_{\veps}^{\ast}$,
  \begin{equation*}
    \det \wt{h}_{\veps}(\xi) \geq a \Lambda_{0, \veps}^{2d}(\xi).
  \end{equation*}
\end{assump}

The Assumptions~\ref{assump:stabatom} and \ref{assump:stabCB} will be
assumed in the sequel without further indication.

\begin{remark}
These assumptions are quite natural and physical. In fact,
Assumption~\ref{assump:stabatom} is just the phonon stability
conditions (for simple Bravais lattice) identified
in~\cite{EMing:2007} represented using the notions of
pseudo-difference operators. Assumption~\ref{assump:stabCB} is the
usual stability condition of a finite element discretization of
continuous problem derived from the Cauchy-Born rule. We note that
as a consequence of these stability assumptions, the continuous
Cauchy-Born elasticity problem is also elliptic, as indicated by
Corollary~\ref{coro:CB} below. From a mathematical point of view,
Assumption~\ref{assump:stabatom} and Assumption~\ref{assump:stabCB}
can be seen as the uniform ellipticity of the difference operator.%
\end{remark}

Next, we prove a lower bound for the symbol $\wt{h}_{\qc}$, which is
crucial for the regularity and stability estimates. Let us recall an
inequality proved by Ky Fan:
\begin{theorem}[Ky Fan's determinant inequality \cite{Fan:1950}]
  \label{thm:KyFan}
  Let $A$, $B$ be positive definite matrices, then for any $\lambda
  \in [0, 1]$,
  \[
    \det( \lambda A + (1 - \lambda) B) \geq
    (\det A)^{\lambda} (\det B)^{1-\lambda}.
  \]
\end{theorem}

\begin{coro}\label{cor:lowbound}
  For any $\veps > 0$, $x \in \Omega_{\veps}$ and any $\xi \in
  \LL_{\veps}^{\ast}$, we have
  \[
    \det \wt{h}_{\qc}(x, \xi) \geq \min(a, a_{\at}) \Lambda_{0, \veps}^{2d}(\xi).
  \]
\end{coro}

\begin{proof}
  This is an immediate corollary of Theorem~\ref{thm:KyFan}. Since for
  any $x$, $\varrho(x) \in [0, 1]$, we have
  \begin{align*}
    \det \wt{h}_{\qc}(x, \xi) & = \det \bigl( (1 - \varrho(x))
    \wt{h}_{\at}(\xi) + \varrho(x) \wt{h}_{\veps}(\xi) \bigr) \\
    & \geq (\det \wt{h}_{\at}(\xi))^{1 - \varrho(x)}
    (\det \wt{h}_{\veps}(\xi))^{\varrho(x)} \\
    & \geq a_{\at}^{1-\varrho(x)} a^{\varrho(x)} \Lambda_{0, \veps}^{2d}(\xi) \\
    & \geq \min(a, a_{\at}) \Lambda_{0, \veps}^{2d}(\xi).
  \end{align*}
\end{proof}

With these preparations, we now establish the regularity estimate of
the quasi-continuum approximation. The regularity of discrete
elliptic systems is understood by a fundamental result of finite
difference approximation by Bube and
Strikwerda~\cite{BubeStrikwerda:1983}. They extended the
regularity estimate of Thom{\'e}e and
Westergren~\cite{ThomeeWestergren:1968} from single elliptic
equation to elliptic systems.

Let us introduce the regular discrete elliptic system
following~\cite{BubeStrikwerda:1983}. The concept is parallel to the
regular continuous elliptic
system~\cite{AgmonDouglisNirenberg:1959}.
\begin{defn}[Regular discrete elliptic system]
  For $i, j = 1, \cdots, d$, let $L_{ij}$ be a difference operator
  with symbol $l_{ij}(x, \xi)$. The system of difference equations
  \begin{equation}\label{eq:ellipsys}
    \sum_{j=1}^d L_{ij} v_j(x) = f_i(x), \qquad i = 1, \cdots, d,
  \end{equation}
  is a \emph{regular discrete elliptic system}, if there are set of
  integers $\{\sigma_i\}_{i=1}^d$ and $\{\tau_j\}_{j=1}^d$ such that
  each $L_{ij}$ is a difference operator of order at most $\sigma_i +
  \tau_j$, and if there are positive constants $C, \xi_0, \veps_0$
  such that
  \begin{equation*}
    \abs{ \det l_{ij}(x, \xi) } \geq C \Lambda_{\veps}^{2p}(\xi)
  \end{equation*}
  for $ 0 < \veps \leq \veps_0$, $\xi \in \LL_{\veps}^{\ast}$, and
  $\max_{1\leq i \leq d} \abs{\xi_i} \geq \xi_0$, where $2p = \sum_{i}
  (\sigma_i + \tau_i)$. We will say that the system
  \eqref{eq:ellipsys} is regular elliptic of order $(\sigma, \tau)$.
\end{defn}

By Corollary~\ref{cor:lowbound}, we immediately have
\begin{prop}\label{prop:elliptic}
  Under Assumptions~\ref{assump:stabatom} and \ref{assump:stabCB},
  the finite difference system
  \begin{equation}\label{eq:ellipsys2}
    \mc{H}_{\qc} v = f
  \end{equation}
  is a regular discrete elliptic system of order $(0,2)$.
\end{prop}

For the regular discrete elliptic system~\eqref{eq:ellipsys2}, we
have the following regularity estimate.
\begin{theorem}\label{thm:regularity}
  Under Assumptions~\ref{assump:stabatom} and \ref{assump:stabCB}, for
  any $v \in H^2_{\veps}(\Omega)$, we have
  \begin{equation}\label{eq:regularity}
    \norm{v}_{\veps,2} \leq
    C (\norm{\mc{H}_{\qc} v}_{\veps,0} + \norm{v}_{\veps,0}).
  \end{equation}
  The constant $C$ is independent of $v$ and $\veps$.
\end{theorem}

\begin{remark}
  Theorem~\ref{thm:regularity} is analogous to the interior regularity
  estimate for elliptic partial differential equations given in
  \cite{AgmonDouglisNirenberg:1964}.  The statement of the theorem is
  just rewriting Theorem 2.1 in \cite{BubeStrikwerda:1983} using the
  current notation. We note that in \cite{BubeStrikwerda:1983}, Bube
  and Strikwerda proved interior regularity estimates, which clearly
  implies the \textit{a priori} estimate for periodic case here.
\end{remark}
\section{Stability}\label{sec:stability}

The main theorem we will prove in this section is the following
stability estimate.
\begin{theorem}[Stability]\label{thm:stability}
  Under Assumptions~\ref{assump:stabatom} and \ref{assump:stabCB}, for
  any $v \in H^2_{\veps}(\Omega)$, we have
  \begin{equation}\label{eq:stability}
    \norm{v}_{\veps,2} \leq C \norm{\mc{H}_{\qc} v}_{\veps,0}.
  \end{equation}
\end{theorem}

Let us make some remarks about the stability result. In general, we do
not know whether a stability estimate like \eqref{eq:stability} is
valid for the force-based quasicontinuum method in general dimension
(see \cites{DobsonLuskinOrtner:2010b, DobsonOrtnerShapeev} for some
study in one dimension). From a pseudo-difference operator point of
view, the continuity in $x$ variable of the symbol of the linearized
operator is crucial for the validity of the strong stability. This is
also the main motivation to use a smooth transition function
$\varrho(x)$ in the current scheme. The strong stability property of
the scheme will facilitate the numerical solution based on iterative
methods.

We also note that the strong stability is also crucial for the
extension of the current scheme to the time-dependent case. It plays
the role of G\r{a}rding inequality.  We will leave this to future
publications.

To obtain the stability estimate from the regularity estimate of
Theorem~\ref{thm:regularity}, we need to eliminate $\norm{v}_{\veps,
  0}$ on the right hand side of \eqref{eq:regularity}. In spatial
dimension one, this can be achieved by the discrete maximum
principle for the finite difference equation. This is however no
longer the case for higher dimensions, as then we are dealing with
an elliptic system. The argument we will use is instead similar in
spirit to the argument used in \cites{AgmonDouglisNirenberg:1959,
Schechter:1959} for passing from regularity estimate to uniqueness
results for elliptic systems.

The difficulty however is that a compactness argument as in
\cite{Schechter:1959} can not apply to the finite difference system,
as we need a uniform estimate for different $\veps$. Therefore,
instead of using the compactness, the proof is based on the uniqueness
of the continuous system from ellipticity, the consistency of the
finite difference schemes to the continuous system, and the regularity
estimate Theorem~\ref{thm:regularity}. We note that a similar approach
was considered by Martin \cite{Martin:94}.

In order to connect the finite difference system with continuous PDE,
we need to extend grid functions on $\Omega_{\veps}$ to continuous
functions defined in $\Omega$. For this purpose, let us define an
interpolation operator $Q_{\veps}$ as follows.\footnote{ Usual linear
  interpolations are not sufficient for our purpose as we need high
  regularity of the interpolated functions.} For any lattice function $u$ on
$\Omega_{\veps}$, we define $Q_{\veps} u \in L^2(\Omega)$ as
\begin{equation}\label{eq:Qvepsdef}
  (Q_{\veps} u)(x) = (2\pi)^{d/2} \sum_{\xi \in \mathbb{L}^{\ast}_{\veps}}
  e^{\I x \cdot \xi} \wh{u}(\xi), \quad x \in \Omega.
\end{equation}
Comparing with \eqref{eq:invdiscF}, we know that $Q_{\veps} u$ agrees
with $u$ on $\Omega_{\veps}$. We have the following 
properties of $Q_{\veps}$.
\begin{lemma} For $k \geq 0$, there exists constants $c_k, C_k > 0$,
  such that for any $u$,
  \begin{equation*}
    c_k \norm{u}_{H^k_{\veps}(\Omega)} \leq
    \norm{Q_{\veps} u}_{H^k(\Omega)} \leq C_k \norm{u}_{H^k_{\veps}(\Omega)}.
  \end{equation*}
\end{lemma}
\begin{proof}
  The conclusion follows immediately from definition
  \eqref{eq:Qvepsdef} and \eqref{eq:symbolcompare}.
\end{proof}

Let $\chi$ be a standard nonnegative cut-off function on $\RR^d$,
which is smooth and compactly supported, with $\norm{\chi}_{L^1} = 1$.
Let $\chi_{\veps}$ be the scaled version
\begin{equation*}
  \chi_{\veps}(x) = \veps^{-(\alpha d)} \chi( \veps^{-\alpha} x),
\end{equation*}
for some $\alpha$ with $0 < \alpha < 1$. The choice of the value of
$\alpha$ will be specified later in the proof of
Proposition~\ref{prop:contcons}.

Define a low-pass filter operator $L_{\veps}$ for $f \in L^2(\Omega)$
using $\wh{\chi_{\veps}}$ as Fourier multiplier:
\begin{equation*}
  \wh{L_{\veps} f}(\xi) = (2\pi)^{d/2} \wh{f}(\xi) \wh{\chi_{\veps}}(\xi)
  = (2\pi)^{d/2} \wh{f}(\xi) \wh{\chi}(\veps^{\alpha} \xi).
\end{equation*}
In real space, $L_{\veps}$ convolves $f$ with $\chi_{\veps}$.  Note
that, using integration by parts, it is easy to see that
\begin{align}
  & \label{eq:chidecay} \abs{\wh{\chi_{\veps}}(\xi)} \leq C_k
  \abs{\veps^{\alpha} \xi}^{-k}, \quad \forall\ k \in \ZZ_+, \\
  & \label{eq:chiunity} (2\pi)^{d/2} \wh{\chi_{\veps}}(0) = 1.
\end{align}
Hence, $L_{\veps}$ is indeed a low-pass filter.
For simplicity of notation, we will denote
\begin{equation*}
  \wb{u}_{\veps} = L_{\veps} Q_{\veps} u_{\veps},
\end{equation*}
for lattice function $u_{\veps}$ on $\Omega_{\veps}$.

We state and prove a consistency result for the linearized operator in
terms of symbols.
\begin{prop}[Consistency of linearized operator]\label{prop:Hcons}
  There exists $\veps_0 > 0$ and $s > 0$ such that for any $\veps \leq
  \veps_0$ and $\xi, \eta \in \LL^{\ast}_{\veps}$, we have
  \begin{equation*}
    \abs{\wh{h}_{\CB}(\xi, \eta) - \wh{h}_{\qc}(\xi, \eta)} \leq C \veps^2
    (\abs{\eta} + 1)^s.
  \end{equation*}
\end{prop}

\begin{proof}
  By definition, for $1 \leq j, k \leq d$,
  \begin{equation*}
    \begin{aligned}
      (\wh{h}_{\qc})_{jk}(\xi, \eta) & = \veps^d (2\pi)^{-d/2} \sum_{x \in
        \Omega_{\veps}} e^{-\I \xi \cdot x} (\wt{h}_{\qc})_{jk}(x, \eta) \\
      & = \veps^d (2\pi)^{-d/2} \sum_{x \in \Omega_{\veps}} e^{-\I
        (\xi+\eta) \cdot x} (\mc{H}_{\qc} (e_k f_{\eta}))_j(x).
    \end{aligned}
  \end{equation*}
  where $f_{\eta}(x) = e^{\I x\cdot \eta}$ for $x \in \Omega$ and
  \begin{equation*}
    \begin{aligned}
      (\wh{h}_{\CB})_{jk}(\xi, \eta) & = (2\pi)^{-d/2} \int_{\Omega} e^{- \I
        \xi\cdot x} \ud x (\wt{h}_{\CB})_{jk}(\eta) \\
      & = \veps^d (2\pi)^{d/2} \sum_{x \in \Omega_{\veps}} e^{-\I \xi
        \cdot x} (\wt{h}_{\CB})_{jk}(\eta) \\
      & = \veps^d (2\pi)^{-d/2} \sum_{x \in \Omega_{\veps}} e^{-\I
        (\xi+\eta) \cdot x} (\mc{H}_{\CB} (e_k f_{\eta}))_j(x),
    \end{aligned}
  \end{equation*}
  where we have used in the fact that $\wt{h}_{\CB}(x, \eta) =
  \wt{h}_{\CB}(\eta)$ due to translational symmetry. Note that we get
  from the second line from the first line in the above equation using
  the fact that $\xi$ takes value in $\LL^{\ast}_{\veps}$, so that the
  integral equals to the sum.

  Hence, taking difference of the above two equations, we obtain the
  bound
  \begin{equation*}
    \abs{\wh{h}_{\qc}(\xi, \eta) - \wh{h}_{\CB}(\xi, \eta) } \leq C
    \sup_{1\leq k \leq d} \norm{\mc{H}_{\qc} (e_k f_{\eta}) - \mc{H}_{\CB}
      (e_k f_{\eta})}_{L^{\infty}_{\veps}}.
  \end{equation*}
  Note that by the definition of linearized operators $\mc{H}_{\qc}$
  and $\mc{H}_{\CB}$, we have
  \begin{equation*}
    \mc{H}_{\qc} (e_k f_{\eta}) - \mc{H}_{\CB} (e_kf_{\eta})
    = \lim_{t \to 0^+} \frac{1}{t} \bigl( \mc{F}_{\qc}[x + t (e_k f_{\eta})] -
    \mc{F}_{\CB}[x + t (e_k f_{\eta})] \bigr).
  \end{equation*}
  Hence,
  \begin{equation*}
    \begin{aligned}
      \norm{\mc{H}_{\qc} (e_kf_{\eta}) - \mc{H}_{\CB}
        (e_kf_{\eta})}_{L^{\infty}_{\veps}} & = \lim_{t \to 0^+}
      \frac{1}{t} \norm{\mc{F}_{\qc}[x + t (e_kf_{\eta})] -
        \mc{F}_{\CB}[x + t (e_kf_{\eta})]}_{L^{\infty}_{\veps}} \\
      & \leq C \veps^2 \norm{e_k f_{\eta}}_{W^{16, \infty}} \leq C
      \veps^2 \norm{e_k f_{\eta}}_{H^{s}} \leq C \veps^2 (1 +
      \abs{\eta})^{s},
    \end{aligned}
  \end{equation*}
  where $s$ is chosen so that the Sobolev inequality
  \begin{equation*}
    \norm{f}_{W^{16, \infty}(\Omega)} \leq C \norm{f}_{H^{s}(\Omega)}
  \end{equation*}
  holds for any $f \in H^{s}(\Omega)$ ($s$ depends on the dimension).
  Here, we have used Corollary~\ref{coro:consqc}, noticing that
  $\norm{t e_k f_{\eta}}_{L^{\infty}}$ is uniformly bounded for $\eta$
  as $t \to 0$. This concludes the proof.
\end{proof}

The proof of Proposition~\ref{prop:Hcons} actually gives for any
  $\veps \leq \veps_0$, $x \in \Omega_{\veps}$ and $\eta \in
  \LL_{\veps}^{\ast}$,
  \begin{equation}\label{eq:Hcons}
    \abs{ \wt{h}_{\qc}(x, \eta) - \wt{h}_{\CB}(\eta) } \leq C \veps^2
    (1 + \abs{\eta})^s.
  \end{equation}
  Combined with Corollary~\ref{cor:lowbound}, we get as a corollary
\begin{coro}\label{coro:CB}
  $\wt{h}_{\CB}(\xi)$ is positive definite and there exists
  $a_{\CB}>0$ such that for any $\xi \in \LL^{\ast}$,
  \begin{equation*}
    \det \wt{h}_{\CB}(\xi) \geq a_{\CB}\Lambda_0^{2d}(\xi).
  \end{equation*}
\end{coro}

\begin{proof}
  Fixed $\xi \in \LL^{\ast}$, take $\veps_1$ sufficiently small, so
  that for $\veps < \veps_1$, $\xi \in \LL_{\veps}^{\ast}$ (it
  suffices to take $\veps_1$ so small that $\veps_1 \xi \in
  \Gamma^{\ast}$). Without loss of generality, we can take $\veps_1$ less than
  $\veps_0$ in Proposition~\ref{prop:Hcons}.

  From the continuous dependence of matrix determinants on matrix
  elements, we get from \eqref{eq:Hcons} that for any $\veps \leq \veps_1$
  sufficiently small, $x \in \Omega_{\veps}$
  \begin{equation*}
    \abs{\det \wt{h}_{hy}(x, \xi) - \det \wt{h}_{\CB}(\xi)}
    \leq C \veps^2 (1 + \abs{\xi})^s.
  \end{equation*}
  Combining the last inequality with Corollary~\ref{cor:lowbound},
  we get the desired estimate by taking $\veps \to 0$.
\end{proof}

With these preparations, let us now state the key proposition
will be used in the proof of Theorem~\ref{thm:stability}.
\begin{prop}\label{prop:contcons}
  For $\{v_{\veps}\}_{\veps > 0}$ that $v_{\veps}\in
  H^2_{\veps}(\Omega)$ and $\norm{v_{\veps}}_{\veps, 2}$ is uniformly
  bounded, we have
  \begin{equation}
    \lim_{\veps \to 0+} \norm{\mc{H}_{\CB} \wb{v}_{\veps}-
      \wb{\mc{H}_{\qc} v_{\veps}}}_{L^2(\Omega)} = 0.
  \end{equation}
\end{prop}

Assume the validity of Proposition~\ref{prop:contcons}, which we will
come back in the end of this section, the proof of
Theorem~\ref{thm:stability} follows a \textit{reductio ad absurdum}.

\begin{proof}[Proof of Theorem~\ref{thm:stability}]

  Suppose \eqref{eq:stability} does not hold, then there is a sequence
  of functions $\{w_k\}$ and $\veps_k > 0$ such that
  \begin{align*}
    & \norm{w_k}_{\veps_k, 2} \to \infty, && \text{as } k \to \infty; \\
    & \norm{\mc{H}_{\qc} w_k}_{\veps_k, 0} \leq c, && \text{for all } k; \\
    & \sum_{x\in \Omega_{\veps_k}} w_k(x) = 0, && \text{for all } k.
  \end{align*}
  Set $v_k = w_k / \norm{w_k}_{\veps_k, 2}$, we then have
  \begin{align}
    \label{eq:useq1} & \norm{v_k}_{\veps_k, 2} = 1 && \text{for all } k;\\
    \label{eq:useq2}
    & \norm{\mc{H}_{\qc} v_k}_{\veps_k, 0} \to 0, && \text{as } k \to \infty; \\
    \label{eq:useq3} & \sum_{x\in \Omega_{\veps_k}} v_k(x) = 0, &&
    \text{for all } k.
  \end{align}
  Since
  \begin{equation*}
    \mc{H}_{\CB} \wb{v}_k = \wb{\mc{H}_{\qc} v_k} +
    ( \mc{H}_{\CB} \wb{v}_k - \wb{\mc{H}_{\qc} v_k}).
  \end{equation*}
  Since $\norm{\mc{H}_{\qc} v_k}_{\veps_k, 0} \to 0$, we have
  \begin{equation*}
    \norm{\wb{\mc{H}_{\qc} v_k}}_{L^2(\Omega)} \to 0, \quad
    \text{as } k \to \infty.
  \end{equation*}
  Moreover, by Proposition~\ref{prop:contcons},
  \begin{equation*}
    \norm{\mc{H}_{\CB} \wb{v}_k -
      \wb{\mc{H}_{\qc} v_k}}_{L^2(\Omega)} \to 0, \quad \text{as } k \to \infty.
  \end{equation*}
  Hence $\norm{\mc{H}_{\CB} \wb{v}_k }_{L^2(\Omega)} \to 0$. Note also
  that the average of $\wb{v}_k$ is zero, since $\wh{\wb{v}_k}(0) =
  0$. By the invertibility of $\mc{H}_{\CB}$ on the subspace
  orthogonal to constant function, $\norm{\wb{v}_k}_{L^2(\Omega)} \to
  0$, as $ k \to \infty$, while $\norm{v_k}_{\veps_k, 2} = 1$. It
  follows then $\norm{v_k}_{\veps_k, 0} \to 0$. Indeed, since
  \begin{equation*}
    \norm{v_k}_{\veps_k, 1} = \sum_{\xi \in \LL^{\ast}_{\veps_k}}
    \Lambda^2_{\veps_k}(\xi) \abs{\wh{v_k}(\xi)}^2 \leq 1,
  \end{equation*}
  for any $\delta > 0$, there exist $\Xi>0$ and $k_1$, such that for
  any $k \geq k_1$,
  \begin{equation}\label{eq:k1}
    \sum_{\xi \in \LL^{\ast}_{\veps_k},\, \abs{\xi} \geq \Xi}
    \abs{\wh{v_k}(\xi)}^2 < \delta / 2.
  \end{equation}
  On the other hand, due to \eqref{eq:chiunity}, there exists $k_2$,
  such that for $k \geq k_2$
  \begin{equation}\label{eq:k2}
    \sum_{\xi \in \LL^{\ast}_{\veps_k}, \, \abs{\xi} < \Xi}
    \bigl\lvert \abs{\wh{v_k}(\xi)^2} - \abs{\wh{\wb{v}_k}(\xi)}^2
    \bigr\rvert \leq \delta / 4.
  \end{equation}
  Moreover, as $\norm{\wb{v}_k}_{L^2} \to 0$, there exists $k_3$, such
  that for $k \geq k_3$,
  \begin{equation}\label{eq:k3}
    \sum_{\xi \in \LL^{\ast}_{\veps_k}, \abs{\xi} < \Xi}
    \abs{\wh{\wb{v}_k}(\xi)}^2 \leq \delta / 4.
  \end{equation}
  Combined \eqref{eq:k1}--\eqref{eq:k3} together, we have for $k \geq
  \max(k_1, k_2, k_3)$,
  \begin{equation*}
    \norm{v_k}_{\veps_k, 0}^2 = \sum_{\xi \in \LL^{\ast}_{\veps_k}}
    \abs{\wh{v_k}}^2\leq \delta.
  \end{equation*}
  Hence, $\lim_{k \to \infty} \norm{v_k}_{\veps_k, 0} = 0$. From
  Theorem~\ref{thm:regularity}, this implies
  \begin{equation*}
    \lim_{k\to \infty} \norm{v_k}_{\veps_k, 2} = 0.
  \end{equation*}
  The contradiction with the choice of $v_k$ proves the Theorem.
\end{proof}

Using perturbation, we may extend the results of
Theorem~\ref{thm:stability} to a deformed state $u$.

\begin{theorem}[Stability]\label{thm:ustability}
  Under Assumption~\ref{assump:stabatom} and \ref{assump:stabCB},
  there exists $\delta > 0$, such that for any $\veps > 0$ and $u$,
  $\norm{u}_{W^{2, \infty}_{\veps}} \leq \delta$ and any $v \in
  H^2_{\veps}(\Omega)$, we have
  \begin{equation}\label{eq:ustability}
    \norm{v}_{\veps,2} \leq C \norm{\mc{H}_{\qc}[u] v}_{\veps,0},
  \end{equation}
  where the constant depends on $\delta$, but is independent of $u$,
  $v$ and $\veps$.
\end{theorem}

\begin{proof}
This theorem follows from a perturbation argument of Theorem~\ref
{thm:stability}. Denote by $v_0$ the solution of
\[
\mc{H}_{\qc}[0] v_0=f.
\]
We immediately have
\[
\mc{H}_{\qc}[0](v-v_0)=\Lr{\mc{H}_{\qc}[0]-\mc{H}_{\qc}[u]}v.
\]
Using Theorem~\ref {thm:stability}, we have
\[
\norm{v-v_0}_{\veps,2}\le
C\norm{\Lr{\mc{H}_{\qc}[0]-\mc{H}_{\qc}[u]}v}_{\veps,0}\le
C \norm{\na u}_{W^{1,\infty}_{\veps}}\norm{v}_{\veps,2}.
\]
By triangular inequality, we have
\begin {align*}
\norm{v}_{\veps,2}&\le\norm{v_0}_{\veps,2}+\norm{v-v_0}_{\veps,2}\\
&\le C\norm{\mc{H}_{\qc}[0] v_0}_{\veps,0}+C \norm{\na
u}_{W^{1,\infty}_{\veps}}\norm{v}_{\veps,2}\\
&=C\norm{\mc{H}_{\qc}[u] v}_{\veps,0}+C\norm{\na
u}_{W^{1,\infty}_{\veps}}\norm{v}_{\veps,2}\\
&\le C\norm{\mc{H}_{\qc}[u] v}_{\veps,0}+C\delta\norm{v}_{\veps,2},
\end {align*}
which gives~\eqref {eq:ustability} by choosing $\delta=1/(2C)$.
\end{proof}

We conclude this section with the proof of
Proposition~\ref{prop:contcons}.
\begin{proof}[Proof of Proposition~\ref{prop:contcons}]
  We work in the Fourier domain. By definition,
  \begin{equation*}
    (\mc{H}_{\CB} \wb{v}_{\veps})(x) =
    \sum_{\xi \in \LL^{\ast}_{\veps}} e^{\I x \cdot \xi} \wt{h}_{\CB}(x, \xi)
    \wh{\chi}(\veps^{\alpha} \xi) \wh{v_{\veps}}(\xi).
  \end{equation*}
  Hence, taking Fourier transform,
  \begin{equation*}
    \begin{aligned}
      \wh{\mc{H}_{\CB} \wb{v}_{\veps}}(\xi) & = (2\pi)^{-d/2}
      \int_{\Omega} e^{-\I \xi \cdot x} \sum_{\eta \in
        \LL^{\ast}_{\veps}} e^{\I x \cdot \eta} \wt{h}_{\CB}(x, \eta)
      \wh{\chi}(\veps^{\alpha} \eta) \wh{v_{\veps}}(\eta) \ud x \\
      & = \sum_{\eta \in \LL^{\ast}_{\veps}} \wh{h}_{\CB}(\xi - \eta, \eta)
      \wh{\chi}(\veps^{\alpha} \eta) \wh{v_{\veps}}(\eta),
    \end{aligned}
  \end{equation*}
  where
  \begin{equation*}
    \wh{h}_{\CB}(\xi, \eta) = (2\pi)^{-d/2} \int_{\Omega} e^{-\I \xi \cdot x}
    \wt{h}_{\CB}(x, \eta) \ud x
  \end{equation*}
  is the Fourier transform of the symbol with respect to $x$.

  On the other hand, for the discrete system, we have
  \begin{equation*}
    \begin{aligned}
      \wh{\wb{\mc{H}_{\qc} v_{\veps}}}(\xi) & = \wh{\chi}(\veps^{\alpha}
      \xi) \veps^d (2\pi)^{-d/2} \sum_{x \in \Omega_{\veps}} e^{-\I
        \xi \cdot x} \sum_{\eta \in \LL^{\ast}_{\veps}} e^{\I x \cdot
        \eta}\wt{h}_{\qc}(x, \eta) \wh{v_{\veps}}(\eta) \\
      & = \wh{\chi}(\veps^{\alpha} \xi) \sum_{\eta \in \LL^{\ast}_{\veps}}
      \wh{h}_{\qc}(\xi - \eta, \eta) \wh{v_{\veps}}(\eta),
    \end{aligned}
  \end{equation*}
  where
  \begin{equation*}
    \wh{h}_{\qc}(\xi, \eta) =
    \veps^d (2\pi)^{-d/2} \sum_{x \in \Omega_{\veps}} e^{-\I \xi \cdot x}
    \wt{h}_{\qc}(x, \eta).
  \end{equation*}

  Let us compare the difference between $\mc{H}_{\CB} \wb{v}_{\veps}$
  and $\wb{\mc{H}_{\qc} v_{\veps}}$. We write
  \begin{equation*}
    \begin{aligned}
      \left\lvert \wh{\mc{H}_{\CB} \wb{v}_{\veps}}(\xi) \right. & -
      \left. \wh{\wb{\mc{H}_{\qc} v_{\veps}}}
        (\xi) \right\rvert  \\
      & = \biggl\lvert \sum_{\eta \in \LL^{\ast}_{\veps}} \Bigl(
      \wh{\chi}(\veps^{\alpha} \eta) \wh{h}_{\CB}(\xi - \eta, \eta) -
      \wh{\chi}(\veps^{\alpha} \xi) \wh{h}_{\qc}(\xi - \eta,
      \eta)\Bigr) \wh{v_{\veps}}(\eta) \biggr\rvert \\
      & \leq \abs{\wh{I_1}(\xi)} + \abs{\wh{I_2}(\xi)},
    \end{aligned}
  \end{equation*}
  where
  \begin{align*}
    & \wh{I_1}(\xi) = \sum_{\eta \in \LL^{\ast}_{\veps}} \bigl(
    \wh{\chi}(\veps^{\alpha} \xi) - \wh{\chi}(\veps^{\alpha} \eta) \bigr)
    \wh{h}_{\CB}(\xi - \eta, \eta) \wh{v_{\veps}}(\eta), \\
    & \wh{I_2}(\xi) = \wh{\chi}(\veps^{\alpha} \xi) \sum_{\eta \in
      \LL^{\ast}_{\veps}} \bigl( \wh{h}_{\CB}(\xi - \eta, \eta) -
    \wh{h}_{\qc}(\xi - \eta, \eta) \bigr) \wh{v_{\veps}}(\eta).
  \end{align*}
  It suffices to prove that $L^2$ norms of $I_1$ and $I_2$ both go to
  zero as $\veps \to 0$. Let us estimate $I_1$ first. By the
  smoothness of $\chi$, we have $\abs{\wh{\chi}(\veps^{\alpha}\xi) -
    \wh{\chi}(\veps^{\alpha}\eta)} \leq C \veps^{\alpha}\abs{\xi - \eta}$,
  hence
  \begin{equation*}
    \abs{\wh{I_1}(\xi)} \leq C \veps^{\alpha} \sum_{\eta\in\LL^{\ast}_{\veps}}
    \abs{\xi - \eta} \bigl\lvert \Lambda^{-2}(\eta) \wh{h}_{\CB}(\xi - \eta, \eta)
    \bigr\rvert \abs{\Lambda^2(\eta) \wh{v_{\veps}}(\eta)}.
  \end{equation*}
  Define $\theta(\xi)$ as
  \begin{equation*}
    \theta(\xi) = \abs{\xi} \sup_{\eta \in \LL^{\ast}} \bigl\lvert
    \Lambda^{-2}(\eta) \wh{h}_{\CB}(\xi, \eta) \bigr\rvert.
  \end{equation*}
  By the smoothness of $\wt{h}_{\CB}(x, \xi)$ with respect to $x$ and the
  fact that $\mc{H}_{\CB}$ is a second order operator, we have
  $\abs{\xi \Lambda^{-2}(\eta) \wh{h}_{\CB}(\xi, \eta)} \leq C
  \abs{\xi}^{-d-1}$ uniformly in $\eta$. Hence, $\theta \in
  l^1(\LL^{\ast})$ as a function of $\xi$. Therefore,
  \begin{equation*}
    \begin{aligned}
      \norm{I_1}_{L^2(\Omega)} = \norm{\wh{I_1}}_{l^2(\LL^{\ast})} & \leq C
      \veps^{\alpha} \norm{\theta}_{l^1(\LL^{\ast})} \biggl( \sum_{\eta
        \in L^{\ast}_{\veps}} \Lambda^4(\eta)
      \abs{\wh{v_{\veps}}(\eta)}^2 \biggr)^{1/2} \\
      & \leq C \veps^{\alpha} \norm{\theta}_{l^1(\LL^{\ast})}
      \norm{Q_{\veps} v_{\veps}}_{H^2(\Omega)} \\
      & \leq C \veps^{\alpha} \norm{\theta}_{l^1(\LL^{\ast})}
      \norm{v_{\veps}}_{H^2_{\veps}(\Omega)},
    \end{aligned}
  \end{equation*}
  where the first inequality results from Young's inequality. This
  proves that $\norm{I_1}_{L^2(\Omega)}$ goes to zero as $\veps \to
  0$.

  Let us consider $I_2$ next. Take $\alpha_1 \in (\alpha, 1)$, we
  break $I_2$ into three parts
  \begin{equation*}
    \wh{I_2}(\xi) = \wh{I_{21}}(\xi) + \wh{I_{22}}(\xi) + \wh{I_{23}}(\xi),
  \end{equation*}
  where
  \begin{align*}
    & \wh{I_{21}}(\xi) = 1_{\abs{\xi} \geq \pi \veps^{-\alpha_1}}
    \wh{\chi}(\veps^{\alpha} \xi) \sum_{\eta \in \LL^{\ast}_{\veps}}
    \bigl( \wh{h}_{\CB}(\xi - \eta, \eta) - \wh{h}_{\qc}(\xi - \eta,
    \eta) \bigr) \wh{v_{\veps}}(\eta), \\
    & \wh{I_{22}}(\xi) = 1_{\abs{\xi} < \pi \veps^{-\alpha_1}}
    \wh{\chi}(\veps^{\alpha} \xi) \sum_{\substack{\eta \in \LL^{\ast}_{\veps}, \\
        \abs{\eta} \geq 2\pi \veps^{-\alpha_1}}} \bigl(
    \wh{h}_{\CB}(\xi - \eta, \eta) - \wh{h}_{\qc}(\xi - \eta,
    \eta) \bigr) \wh{v_{\veps}}(\eta), \\
    & \wh{I_{23}}(\xi) = 1_{\abs{\xi} < \pi \veps^{-\alpha_1}}
    \wh{\chi}(\veps^{\alpha} \xi) \sum_{\substack{\eta \in \LL^{\ast}_{\veps}, \\
        \abs{\eta} < 2\pi \veps^{-\alpha_1}}} \bigl( \wh{h}_{\CB}(\xi -
    \eta, \eta) - \wh{h}_{\qc}(\xi - \eta, \eta) \bigr)
    \wh{v_{\veps}}(\eta).
  \end{align*}
  We will control each term: $I_{21}$ is small due to the decay
  property of $\wh{\chi}$; $I_{22}$ is small since $\xi$ and $\eta$ is
  well separated; $I_{23}$ is small due to consistency.
  \begin{itemize}
  \item[$I_{21}$:] Define $w$ given by
    \begin{equation*}
      \wh{w}(\xi) = \sum_{\eta\in\LL^{\ast}_{\veps}}
      \bigl(\wh{h}_{\CB}(\xi-\eta, \eta) -  \wh{h}_{\qc}(\xi-\eta, \eta)
      \bigr) \wh{v_{\veps}}(\eta).
    \end{equation*}
    We observe that $\wh{w}(\xi)$ is the Fourier transform of
    \begin{equation*}
      w(x) = (\mc{H}_{\CB} Q_{\veps} v_{\veps})(x) - (Q_{\veps}(\mc{H}_{\qc}
      v_{\veps}))(x).
    \end{equation*}
    Hence, $\norm{w}_{L^2(\Omega)} \leq C\norm{v_{\veps}}_{\veps,
      2}$. By \eqref{eq:chidecay}, we have
    \begin{equation*}
      \abs{\wh{\chi}(\veps^{\alpha} \xi)} \leq C_k \veps^{k(\alpha_1 - \alpha)},
      \quad \forall\, \abs{\xi} \geq \pi \veps^{-\alpha_1},
    \end{equation*}
    for any positive integer $k$. Therefore, we conclude that
    $\norm{I_{21}}_{L^2(\Omega)} \to 0$ as $\wh{I_{21}}(\xi) =
    1_{\abs{\xi} \geq \pi \veps^{-\alpha_1}} \wh{\chi}(\veps^{\alpha}\xi)
    \wh{w}(\xi)$.

  \item[$I_{22}$:] We have
    \begin{multline}\label{eq:I22}
      \abs{\wh{I_{22}}(\xi)} \leq C \sum_{\eta \in \LL^{\ast}_{\veps}
      } \abs{\Lambda^{-2}(\eta) \wh{h}_{\CB}(\xi - \eta, \eta)}
      \abs{\Lambda^2(\eta) \wh{v_{\veps}}(\eta)} 1_{\abs{\xi - \eta} >
        \pi\veps^{-\alpha_1}} \\
      + C \sum_{\eta \in \LL^{\ast}_{\veps} } \abs{\Lambda^{-2}(\veps,
        \eta) \wh{h}_{\qc}(\xi - \eta, \eta)}
      \abs{\Lambda^2_{\veps}(\eta) \wh{v_{\veps}}(\eta)} 1_{\abs{\xi -
          \eta} > \pi\veps^{-\alpha_1}}.
    \end{multline}
    The argument for the two terms are analogous, and let us focus on
    the first term. Consider $\varphi(\xi)$ given by
    \begin{equation*}
      \varphi(\xi) = \sup_{\eta \in\LL^{\ast}}
      \abs{\Lambda^{-2}(\eta) \wh{h}_{\CB}(\xi, \eta)}.
    \end{equation*}
    Since $\wt{h}_{\CB}(x, \eta)$ is smooth with respect to $x$ and
    $\mc{H}_{\CB}$ is a second-order operator, we have $\varphi \in
    l^1(\LL^{\ast})$ as a function of $\xi$.  Hence
    \begin{equation*}
      \lim_{\veps\to 0}
      \norm{\varphi(\xi) 1_{\abs{\xi} > \pi \veps^{-\alpha_1}}}_{l^1(\LL^{\ast})} = 0.
    \end{equation*}
    Therefore, using Young's inequality, the first term on the right
    hand side of \eqref{eq:I22} is bounded by $ C \norm{\varphi(\xi)
      1_{\abs{\xi} > \pi \veps^{-\alpha_1}}}_{l^1(\LL^{\ast})}
    \norm{Q_{\veps} v_{\veps}}_{H^2(\Omega)}$, which goes to zero as
    $\veps \to 0$. Hence, $I_{22}$ goes to zero in $L^2$ norm.
  \item[$I_{23}$:] From Proposition~\ref{prop:Hcons}, we have
    \begin{equation*}
      \abs{\wh{h}_{\CB}(\xi, \eta) - \wh{h}_{\qc}(\xi, \eta)} \leq C \veps^2
      (\abs{\eta} + 1)^t
    \end{equation*}
    for some $s \geq 0$. As $\abs{\eta} < 2\pi \veps^{-\alpha_1}$, we have
    \begin{equation*}
      \abs{\wh{h}_{\CB}(\xi, \eta) - \wh{h}_{\qc}(\xi, \eta)} \leq C
      \veps^{(2-s\alpha_1)}.
    \end{equation*}
    Therefore,
    \begin{equation*}
      \begin{aligned}
        \sum_{\xi\in\LL^{\ast}} \abs{\wh{I_{23}}(\xi)}^2
        & \leq C \sum_{\substack{\xi \in \LL^{\ast}, \\
            \abs{\xi} < \pi \veps^{-\alpha_1}}}
        \biggl( \sum_{\substack{\eta \in \LL^{\ast}_{\veps}, \\
            \abs{\eta} < 2\pi \veps^{-\alpha_1}}} \bigl(
        \wh{h}_{\CB}(\xi - \eta, \eta) - \wh{h}_{\qc}(\xi - \eta,
        \eta) \bigr) \wh{v_{\veps}}(\eta) \biggr)^2 \\
        & \leq C \sum_{\substack{\eta \in \LL^{\ast}_{\veps}, \\
            \abs{\eta} < 2\pi \veps^{-\alpha_1}}}
        \abs{\wh{v_{\veps}}(\eta)}^2
        \sum_{\substack{\xi \in \LL^{\ast}, \\
            \abs{\xi} < \pi \veps^{-\alpha_1}}} \bigl \lvert
        \wh{h}_{\CB}(\xi - \eta, \eta) - \wh{h}_{\qc}(\xi - \eta,
        \eta) \bigr\rvert^2 \\
        & \leq C \veps^{4 - (2s + d) \alpha_1}
        \sum_{\substack{\eta \in \LL^{\ast}_{\veps}, \\
            \abs{\eta} < 2\pi \veps^{-\alpha_1}}}
        \abs{\wh{v_{\veps}}(\eta)}^2.
      \end{aligned}
    \end{equation*}
    Hence, by choosing $\alpha_1$ (and also $\alpha$) sufficiently
    small that $ \alpha_1 < 4 / (2s + d)$, we have
    $\norm{I_{23}}_{L^2} \to 0$ as $\veps \to 0$.
  \end{itemize}

  Therefore, to sum up, we have proved both $\norm{I_1}_{L^2(\Omega)}$
  and $\norm{I_2}_{L^2(\Omega)}$ go to zero as $\veps \to 0$. The
  proposition is proved.

\end{proof}
\section{Convergence of the force-based hybrid method}
\label{sec:convergence}

With the consistency and stability results prepared in the last three
sections, we are now ready to prove the main result
Theorem~\ref{thm:main}.  The proof follows the spirit of Strang's
convergence proof of nonlinear finite difference schemes
\cite{Strang:1964}.

As a direct consequence of Corollary~\ref{coro:consqc}, we have the
following
\begin{coro}[Higher order expansion]\label{coro:highorder}
  Under the same assumptions of Theorem~\ref{thm:main}, there exist
  positive constants $\delta$ and $M$, so that for any $p > d$ and
  $ f \in W^{15, p}(\Omega) \cap W^{1, p}_{\sharp}(\Omega) $ with
  $\norm{f}_{W^{15, p}} \leq \delta$,
  denote $\wt{y}=x+u(x)$ with $u$ the solution of the Cauchy-Born
  elasticity problem~\eqref{cb:eq}, we then have
  \[
    \norm{\mc{F}_{\qc}[\wt{y}] - f}_{L^{\infty}_{\veps}} \leq M \veps^2.
  \]
\end{coro}
\begin{remark}
  Using the remark under Lemma~\ref{lem:consCB}, the regularity
  assumption of $f$ can be relaxed to $W^{5, p}(\Omega)$ with $p>d$.
\end{remark}

\begin{proof}[Proof of Theorem~\ref{thm:main}]

We take $\wt{y}$ be that given by Corollary~\ref{coro:highorder}. It
is easy to see
\[
  \int_0^1\mc{H}_{\qc}[ty+(1-t)\wt{y}](x)\dt\cdot(y-\wt{y})
  =\mc{F}_{\qc}[y]-\mc{F}_{\qc}[\wt{y}].
\]
Hence $y$ is the solution of~\eqref {sqc:eq} if and only if
\[
  \int_0^1\mc{H}_{\qc}[ty+(1-t)\wt{y}](x)\dt\cdot(y-\wt{y})
  =f-\mc{F}_{\qc}[\wt{y}].
\]
For any $\kappa\in (3/2, 2)$, we define
\[
  B=\set{y\in X_{\veps}}{\norm{y-\wt{y}}_{\veps,2}\le\veps^{\kappa}}.
\]
We define a map $T: B \to B$ as follows: for any $y\in B$, let
$T(y)$ be the solution of the linear system
\begin {equation}\label {linear}
    \int_0^1\mc{H}_{\qc}[ty+(1-t)\wt{y}](x)\dt\cdot\Lr{T(y)-\wt{y}}
    =f-\mc{F}_{\qc}[\wt{y}].
\end {equation}
We first show that $T$ is well defined. Since
\[
  \norm{ty+(1-t)\wt{y}-\wt{y}}_{\veps,2}\le
  t\norm{y-\wt{y}}_{\veps,2}\le\veps^{\kappa},
\]
which gives that for sufficiently small $\veps$ and $d\le 3$, there
holds
\[
  \norm{ty+(1-t)\wt{y}-\wt{y}}_{W_{\veps}^{2,\infty}}\le\veps^{\kappa-d/2}
  <\delta,
\]
  where the constant $\delta$ appears in Theorem~\ref{thm:ustability}.
  It follows from Theorem~\ref{thm:ustability} that the problem~\eqref
  {linear} is solvable and
  \begin {equation}\label{err:eq}
    \begin{aligned}
      \norm{T(y)-\wt{y}}_{\veps,2}&\leq
      C\norm{f-\mc{F}_{\qc}[\wt{y}]}_{\veps,0}\\
      &\le C\norm{f-\mc{F}_{\at}[\wt{y}]}_{\veps,0}
      +C\norm{\mc{F}_{\at}[\wt{y}]-\mc{F}_{\qc}[\wt{y}]}_{\veps,0}\\
      &\le C\veps^2,
    \end{aligned}
  \end{equation}
  where we have used Corollary~\ref{coro:highorder}. For sufficiently
  small $\veps$, we have
  \[
  \norm{T(y)-\wt{y}}_{\veps,2}\le\veps^\kappa.
  \]
  Therefore, $T(y)\in B$ and $T$ is well-defined, which in turn
  implies $T(B)\subset B$ for sufficiently small $\veps$. Now the
  existence of $y$ follows from the Brouwer fixed point theorem. The
  solution $y$ is locally unique since the Hessian at $y$ is
  nondegenerate. Let us denote the solution as $y_{\qc}$, we then have
  from \eqref{err:eq} that
  \begin{equation}\label{err:eq2}
    \norm{\wt{y}-y_{\qc}}_{\veps,2}\le C\veps^2.
  \end{equation}
  Proceeding along the same line that leads to~\eqref {err:eq} and
  using Lemma~\ref{lem:consCB}, we get
  \begin {equation}\label{err:highorder}
    \norm{\wt{y}-y_{\at}}_{\veps,2}\le C\veps^2.
  \end {equation}
  Finally, we conclude that $y_{\qc}$ satisfies~\eqref {eq:main} by
  combining~\eqref{err:eq2} and \eqref{err:highorder}.
\end{proof}

\bibliographystyle{amsalpha}
\bibliography{quasicont}

\end{document}